\documentclass[12pt]{amsart}

\usepackage{amssymb,amsbsy}
\usepackage{latexsym}
\usepackage{verbatim,color}
\usepackage{fullpage}

\input {cyracc.def}
\font\ququ=cmr10 scaled \magstep1
 \font\tencyr=wncyr10 
\font\tencyb=wncyb10 
\font\tencyi=wncyi10 
\font\tencysc=wncysc10 
\def\rus{\tencyr\cyracc}
\def\rusb{\tencyb\cyracc}
\def\rusi{\tencyi\cyracc}
\def\rusc{\tencysc\cyracc}

\renewenvironment{proof}
{\noindent {\sl Proof.}\quad }{\hfill
$\square$ \vskip1.1ex\noindent }

\newenvironment{proof*}
{\noindent {\sl Proof.}\quad }{\hfill
$\square$}

\renewcommand{\theequation}{\thesection .\arabic{equation}}
\renewcommand{\thesubsubsection}{\theequation .\arabic{subsubsection}}

\catcode`\@=\active
\catcode`\@=11
\def\@eqnnum{\hbox to
.01pt{}\rlap{\hskip-\displaywidth(\mathbf{\theequation})}}
\catcode`\@=12

\newenvironment{s}[1]
{ \vskip1.2ex \refstepcounter{equation}
\noindent {\bf \theequation\enspace #1.} \begin{sl}}{\end{sl}
\vskip1.1ex\noindent }

\newenvironment{rem}[1]
{ \vskip1.2ex \refstepcounter{equation}
\noindent {\bf \theequation\enspace {#1}.} }{ \vskip1.1ex\noindent }

\newcommand {\ah}{{\frak a}}

\newcommand {\ce}{{\frak c}}

\newcommand {\g}{{\frak g}}
\newcommand {\h}{{\frak h}}
\newcommand {\ka}{{\frak k}}

\newcommand {\ma}{{\frak m}}

\newcommand {\q}{{\frak q}}

\newcommand {\es}{{\frak s}}
\newcommand {\te}{{\frak t}}

\newcommand {\fF}{{\frak F}}

\newcommand {\z}{{\frak z}}
\newcommand {\gln}{{\frak gl}_n}
\newcommand {\sln}{{\frak sl}_n}

\newcommand {\sono}{{\frak so}_{2n+1}}

\newcommand {\son}{{\frak so}_{n}}

\newcommand {\gN}{{\goth N}}

\newcommand {\esi}{\varepsilon}

\newcommand {\ap}{\alpha}

\newcommand {\N}{{\mathcal N}}
\newcommand {\co}{{\mathcal O}}

\newcommand {\BZ}{{\mathbb Z}}

\newcommand {\BR}{{\mathbb R}}

\newcommand {\md}{/\!\!/}

\newcommand {\ad}{{\mathrm{ad\,}}}
\newcommand {\ads}{{\mathrm{ad}^*}}
\newcommand {\Ad}{{\mathrm{Ad\,}}}
\newcommand {\Ads}{{\mathrm{Ad}^*}}
\newcommand {\Ann}{{\mathrm{Ann\,}}}
\newcommand {\bideg}{{\mathrm{bideg\,}}}
\newcommand {\codim}{{\mathrm{codim\,}}}

\newcommand {\gr}{{\mathrm{gr\,}}}

\newcommand {\ind}{{\mathrm{ind\,}}}
\newcommand {\Lie}{{\mathrm{Lie\,}}}

\newcommand {\Ima}{{\mathrm{Im\,}}}
\newcommand {\Mor}{\operatorname{Mor}}

\newcommand {\rk}{{\mathrm{rk\,}}}
\newcommand {\spe}{{\mathrm{Spec\,}}}

\newcommand {\tr}{{\mathrm{tr\,}}}
\newcommand {\trdeg}{{\mathrm{trdeg\,}}}

\newcommand {\tri}{{\frak sl}_2}
\newcommand {\GR}[2]{{\textrm{{\bf #1}}}_{#2}}

\newcommand {\ov}{\overline}
\newcommand {\un}{\underline}

\newcommand {\beq}{\begin{equation}}
\newcommand {\eeq}{\end{equation}}

\renewcommand{\le}{\leqslant}
\renewcommand{\ge}{\geqslant}

\newfam\Bbbfam\newfam\eufam\newfam\eusfam%
\font\Bbbfont=msbm10 scaled 1200%
\font\bbbfont=msbm10 scaled 1000%
\font\olala=msam10 scaled 1200%
\font\frak=eufm10 scaled 1400%
\font\Bbbsmallfont=msbm8%
\textfont\Bbbfam=\Bbbfont\scriptfont\Bbbfam=\Bbbsmallfont
\font\euzw=eufm10 scaled 1200%
\font\euac=eufm7 scaled 1200%
\font\euacc=eufm7 scaled 1000%
\textfont\eufam=\euzw\scriptfont\eufam=\euac%
\scriptscriptfont\eufam=\euacc%
\font\euszw=eusm10 scaled 1200%
\font\eusac=eusm7 scaled 1200%
\font\eusacc=eusm7 scaled 1000%
\textfont\eusfam=\euszw\scriptfont\eusfam=\eusac%
\scriptscriptfont\eusfam=\eusacc%
\def\frak{\fam\eufam}%
\def\goth{\fam\eusfam}%

\def\square{\hbox {\olala\char"03}}
\def\bbk{\hbox {\Bbbfont\char'174}}
\def\bkk{\hbox {\bbbfont\char'174}}

\begin{document}
\setlength{\parskip}{2pt plus 4pt minus 0pt}
\hfill {\scriptsize \today}
\vskip1ex
\vskip1ex

\title[On coadjoint representations]{On the coadjoint representation of
$\mathbb Z_2$-contractions of reductive Lie algebras}
\author{Dmitri I. Panyushev}
\thanks{This research was supported in part by 
RFBI Grants 05--01--00988 and 06--01--72550}
\maketitle
\vskip-1.5ex
\begin{center}
{\footnotesize
{\it Independent University of Moscow,
Bol'shoi Vlasevskii per. 11 \\
119002 Moscow, \quad Russia \\ e-mail}: {\tt panyush@mccme.ru }\\
}
\end{center}
\tableofcontents
\section*{Introduction}

\noindent
The ground field $\bbk$ is algebraically closed and of characteristic zero.
Let $\g$ be a reductive algebraic Lie algebra. Classical results of Kostant
\cite{ko63} give a fairly complete invariant-theoretic picture of the (co)adjoint
representation of $\g$.
Let $\sigma\in\text{Aut}(\g)$ be an involution and 
$\g=\g_0\oplus\g_1$ the corresponding $\BZ_2$-grading.
Associated to this decomposition, 
there is a non-reductive Lie algebra $\ka=\g_0\ltimes\g_1$, the semi-direct product
of the Lie algebra $\g_0$ and $\g_0$-module $\g_1$. Let $K$ denote a connected
group with Lie algebra $\ka$. 
A remarkable property of the  Lie algebra contraction $\g \leadsto \ka$ is  that
it preserves the transcendence degree of the algebras
of invariants for both adjoint and coadjoint representations of $\ka$; i.e.,
$\trdeg \bbk[\ka]^K =\trdeg \bbk[\ka^*]^K=\rk\g$. The latter equality also shows that
$\ind\ka=\rk\g$.
In \cite{p05}, we proved that many good properties of Kostant's picture for
$(\g,\ad)$ carry over to $(\ka,\ad)$. In particular, 
$\bbk[\ka]^K$ is a polynomial algebra and the quotient mapping 
$\pi_{\ka}: \ka\to \ka\md K=\spe (\bbk[\ka]^K)$ is equidimensional.
The goal of this article is to study 
the invariants of $(\ka,\ads)$. Motivated by several examples,
we come up with the following

\begin{s}{Conjecture}   \label{conj1}
The algebra of invariants of $(\ka,\ads)$ is polynomial and 
the quotient mapping 
$\pi_{\ka^*}: \ka^*\to \ka^*\md K=\spe (\bbk[\ka^*]^K)$ is equidimensional.
\end{s}%
If $\es$ is reductive, $\g=\es\dotplus\es$, and $\sigma\in\text{Aut}(\g)$ is the permutation, 
then $\ka=\es\ltimes\es$ is a so-called {\it Takiff\/} Lie algebra. Here 
$\ad\simeq\ads$ and the validity of the conjecture
follows from results of Takiff~\cite{takiff} and Geoffriau~\cite{geof} (see also \cite{p05}).
Therefore one can concentrate on the case in which $\g$ is simple, where
the adjoint and coadjoint representations of $\ka$ are different. 
It is not hard to prove that if $\g_1$ contains a Cartan subalgebra of $\g$
(the "maximal rank" case), then the conjecture is true. More generally, 
we prove Conjecture~\ref{conj1} (and some stronger assertions) for the $\N$-{\it regular} 
$\BZ_2$-gradings, i.e., if
$\g_1$ contains a regular nilpotent element of $\g$ (see Section~\ref{sect:N-reg}). 
There are also several cases, where we can
prove only "half" of the conjecture, i.e., the fact that $\bbk[\ka^*]^K$ is polynomial
(Section~\ref{good_systems}).
Our proofs of polynomiality in Section~\ref{good_systems} make use of some general results on coadjoint representations. We show that $(\ka,\ads)$ has
a so-called {\it codim--2 property\/}, i.e., the set of non-regular elements of $\ka^*$
is of codimension $\ge 2$ (Theorem~\ref{thm:codim2}). 
This property implies, in turn, that if $l=\ind\ka$ and $F_1,\dots,F_l\in \bbk[\ka^*]^K$ are 
homogeneous and algebraically independent, then 
\begin{equation} \label{intro-ner}
\sum_{i=1}^l\deg F_i \ge (\dim\ka+l)/2\  .
\end{equation}
Furthermore, if the equality holds, then
$F_1,\dots,F_l$ freely generate $\bbk[\ka^*]^K$ (see Theorem~\ref{ner-vo}). 
That is, the polynomiality follows if one could find algebraically independent $K$-invariants 
with "sufficiently small" degrees. To this end, we use the method of 
$\BZ_2$-{\it degeneration\/}
of $G$-invariants in $\bbk[\g]$. Namely, the decomposition
$\g=\g_0\oplus\g_1$ determines the natural bi-grading of $\bbk[\g]$.
For a homogeneous  $f\in\bbk[\g]^G$, let $f^\bullet$ be the bi-homogeneous component
of $f$ of highest degree with respect to $\g_1$. Then regarded as function on $\ka^*$,
$f^\bullet$ is  $K$-invariant (Proposition~\ref{z2-contra}).
Notice that $\deg f^\bullet=\deg f$. It is also known that if 
$f_1,\dots,f_l\in \bbk[\g]^G$ are basic invariants, then 
$\sum_{i=1}^l\deg f_i = (\dim\g+l)/2$. Hence it suffices to find
a set of basic $G$-invariants $f_1,\dots,f_l$ such that 
$f^\bullet_1,\dots,f^\bullet_l$ are algebraically independent. 
In this situation, we say that $f_1,\dots,f_l$ form a {\it good generating
system\/} for $(\g,\g_0)$. Then
the functions $\{ f^\bullet \mid f\in \bbk[\g]^G\}$ form the whole algebra
$\bbk[\ka^*]^K$. However, this is not always the case
 (see Remark~\ref{bad}). 
Therefore the proof of polynomiality
for some symmetric pairs requires different ideas.
\\[.6ex]
Inequality~\eqref{intro-ner} holds for any Lie algebra with the codim--2 property.
But $\bbk[\ka^*]^K$ is bi-graded (Theorem~\ref{semidir_gen}), and we also prove 
a bi-graded refinement of that inequality (see Theorem~\ref{bideg}). In the last section,
we gather the available inforamtion on the bi-degrees of basic invariants for 
$\bbk[\ka]^K$ and $\bbk[\ka^*]^K$.

All Lie algebras are assumed to be algebraic. Algebraic groups are denoted by 
capital Latin letters. Corresponding Lie algebras are denoted by the lowercase Gothic
letters. 
\\[.6ex]
If an algebraic group $Q$ acts on an irreducible affine variety $X$, then $\bbk[X]^Q$ 
is the algebra of $Q$-invariant regular functions on $X$ and $\bbk(X)^Q$
is the field of $Q$-invariant rational functions. If $\bbk[X]^Q$
is finitely generated, then $X\md Q:=\spe \bbk[X]^Q$, and
the {\it quotient morphism\/} $\pi_X: X\to X\md Q$ is the mapping associated with
the embedding $\bbk[X]^Q \hookrightarrow \bbk[X]$. If $\bbk[X]^Q$ is graded polynomial, 
then the elements of any set of algebraically independent homogeneous generators 
will be referred to as {\it basic invariants\/}.
Occasionally, we write $\textit{Inv}(\q,\ad)$ and $\textit{Inv}(\q,\ads)$ for the algebras 
of invariants of the adjoint and coadjoint representations of $\q=\Lie Q$, 
respectively.  If $V$ is a $Q$-module, then $\q_v$ is the stabiliser of $v\in V$ in $\q$.
For the adjoint representation of $\q$, the stabiliser of $x\in \q$ is also denoted by 
$\z_\q(x)$.
A direct sum of Lie algebras is denoted by '$\dotplus$'.

Given an irreducible variety $Y$,  an open subset
$\Omega\subset Y$ is said to be {\it big\/} if $Y\setminus\Omega$ contains no divisors.

$[n]=\{1,2,\dots,n\}$; \ $\lfloor x\rfloor$ is the least integer not exceeding $x$.

\noindent
{\bf Acknowledgements.} {\small 
I wish to thank Sasha Premet and Oksana Yakimova for sharing some important insights
and enlightening discussions on coadjoint representations.
Part of this work was done during my 
stay at the Max-Planck-Institut f\"ur Mathematik (Bonn). 
I am grateful to this institution for the warm hospitality and support.}


\section{The codim-2 property for coadjoint representations}
\label{sect:gen_co}
\setcounter{equation}{0}

\noindent
Let $Q$ be a connected algebraic group with Lie algebra $\q$.
Let $\q^*_{reg}$ denote the set of all $Q$-{\it regular\/} elements of $\q^*$. That is,
\[
   \q^*_{reg}=\{\xi\in\q^*\mid \dim Q{\cdot}\xi\ge \dim Q{\cdot}\eta 
   \text{ for all } \eta\in\q^*\} \ .
\]
As is well-known, $\q^*_{reg}$ is a dense open subset of $\q^*$.

\begin{rem}{Definition}   \label{def:codim2}
We say that the coadjoint representation of $\q$ has the {\it codim--2 property\/} if
$\codim (\q^*\setminus \q^*_{reg})\ge 2$, i.e., $\q^*_{reg}$ is big.
\end{rem}%
{\bf Example.} If $\g$ is reductive, then $\ad\simeq\ads$ and 
$\codim (\g\setminus \g_{reg})=3$. Hence the coadjoint representation of a 
reductive Lie algebra has the codim--2 property.

If $\xi\in \q^*_{reg}$, then $\dim\q_\xi$ is called the {\it index\/} of $\q$, denoted $\ind\q$.
Recall that each orbit $Q{\cdot}\xi$ is a symplectic variety and hence $\dim Q{\cdot}\xi$ is even.
By Rosenlicht's theorem, $\trdeg\bbk(\q^*)^Q=\ind\q$. It follows that
if $f_1,\dots,f_r\in \bbk[\q^*]^Q$ are algebraically independent, then $r\le \ind\q$.

Importance of the codim--2 property is explained by the following result,
which makes use of some ideas of \cite[Theorem\,3.1]{or} (cf. also \cite[Theorem\,1.2]{prya}).

\begin{s}{Theorem}  \label{ner-vo}
Suppose that $(\q,\ads)$ has the codim--2 property
and $\trdeg \bbk[\q^*]^Q=\ind\q$. Set  $l=\ind\q$. 
Let $f_1,\dots,f_l\in \bbk[\q^*]^Q$ be arbitrary
homogeneous algebraically independent polynomials.
Then 
\begin{itemize}
\item[\sf (i)] \ $\sum_{i=1}^l\deg f_i \ge (\dim\q+\ind\q)/2$;
\item[\sf (ii)] \  If\/  $\sum_{i=1}^l\deg f_i = (\dim\q+\ind\q)/2$, then
$\bbk[\q^*]^Q$ is freely generated by $f_1,\dots,f_l$ and
$\xi\in\q^*_{reg}$ if and only if $(\textsl{d}f_1)_\xi,\dots,(\textsl{d}f_l)_\xi$ are linearly
independent.
\end{itemize}
\end{s}\begin{proof} 
Recall that $\mathcal S(\q)=\bbk[\q^*]$ is a Poisson algebra, and the symplectic leaves
in $\q^*$ are precisely the coadjoint orbits of $Q$. Let $\{\ ,\ \}$ denote the Poisson
bracket in $\bbk[\q^*]$. Then $\bbk[\q^*]^Q$ is the centre of $(\bbk[\q^*], \{\ ,\ \})$.

Let $\pi$ denote the Poisson tensor (bi-vector) on $\q^*$.
If $T(\q^*)$ is the tangent bundle of $\q^*$, then
$\pi$ is a section of $\wedge^{2}T(\q^*)$. By definition, if 
$f_1,f_2\in \mathcal S(\q)$, then $\pi(\textsl{d}f_1,\textsl{d}f_2)=\{f_1,f_2\}$.
In particular, if $x,y\in\q$, then $\pi(\textsl{d}x,\textsl{d}y)=[x,y]$.
We regard $\pi$ as an element of
the graded skew-symmetric algebra of polynomial vector fields on $\q^*$. 
Set $n=\dim\q$ and $l=\ind\q$.
Let $\rk \pi_\xi$ denote the rank of the bi-vector $\pi$ at $\xi\in \q^*$. 
It is easily seen that $\rk\pi_\xi=\dim Q{\cdot}\xi$.
Therefore 
\[
   \{\xi\in\q^*\mid \rk\pi_\xi<n-l\}=\q^*\setminus\q^*_{reg} \ .
\]
It follows from the definition of index that 
\[
\mathfrak V_1:=\underbrace{\pi\wedge\pi\wedge\cdots\wedge \pi}_{{(n-l)/2}}
\]
is the maximal nonzero exterior power of
$\pi$. It is an $(n-l)$-vector field on $\q^*$ of degree $(n-l)/2$ and 
$(\mathfrak V_1)_\xi=0$  if and only if $\xi\not\in\q^*_{reg}$.

On the other hand, given algebraically independent polynomials
$f_1,\dots,f_l\in \bbk[\q^*]^Q$, we
get the nonzero differential $l$-form, $\textsl{d}f_1\wedge\ldots\wedge \textsl{d}f_l$,
on $\q^*$. Let $x_1,\dots,x_n$ be a basis for $\q$.
For an $l$-form $\fF$ on $\q^*$, let $\fF^\diamond$ denote the
$(n-l)$-vector field defined by the formula
\[
   \fF^\diamond(\Psi_1,\ldots,\Psi_{n-l})= \frac{\fF\wedge \Psi_1\wedge\ldots 
\wedge\Psi_{n-l}}{\textsl{d}x_1\wedge \textsl{d}x_2\wedge\ldots\wedge \textsl{d}x_n}\ 
\]
for arbitrary differential 1-forms $\Psi_i$. 
It is easily seen that the operation `$\diamond$' does not affect the degree. That is, if
$\fF$ is homogeneous, then so is $\fF^\diamond$, and $\deg \fF=\deg\fF^\diamond$.
Set $\mathfrak V_2:=(\textsl{d}f_1\wedge\ldots\wedge \textsl{d}f_l)^\diamond$.
Thus, both $\mathfrak V_1$ and $\mathfrak V_2$ are nonzero sections of $\wedge^{n-l}T(\q^*)$.
Note that $\deg\mathfrak V_2=\sum_i(\deg f_i-1)$.

For a vector field $\frak v$, let $\imath_{\frak v}$ denote the contraction of 
a section of  $\wedge^{j}T(\q^*)$ with respect to ${\frak v}$. 
Then 
\begin{equation}  \label{eq:V2}
\imath_{\textsl{d}f_j}\mathfrak V_2=
(\textsl{d}f_1\wedge\ldots\wedge \textsl{d}f_l\wedge \textsl{d}f_j)^\diamond=0
 \text{ for all } j\in\{1,\dots,l\} \ .
\end{equation} 
Since each $f_j$ is a central element of the Poisson algebra $\bbk[\q^*]$, we have
$\imath_{\textsl{d}f_j}\pi=0$. It follows that 

\begin{equation}  \label{eq:V1}
\imath_{\textsl{d}f_j}\mathfrak V_1=0  
\text{ for all } j\in\{1,\dots,l\} \ .
\end{equation}
For $\xi\in\q^*$, let $V_\xi\subset T_\xi(\q^*)\simeq \q^*$ denote
the annihilator of the $\bbk$-linear span of 
$\{(\textsl{d}f_i)_\xi\mid i=1,\ldots, l\}$. Consider the open non-empty subset 
$S:=\{\eta\in\q^*\mid (\textsl{d}f_1)_\eta\wedge\ldots\wedge (\textsl{d}f_l)_\eta\ne 0\}$. 
If $\xi\in S$, then $\dim V_\xi=n-l$. 
Let $t\in\wedge^{n-l} T_\xi (\q^*)$ be an
$(n-l)$-vector such that $\imath_{(\textsl{d}f_i)_\xi}t=0$, $i=1,\dots,l$.
Using the undergraduate linear algebra, one readily shows that
$t\in\wedge^{n-l}V_\xi$. Applying this to Eq.~\eqref{eq:V2} and \eqref{eq:V1}, we see
that, for each $\xi\in S$, $(\mathfrak V_1)_\xi$ and $(\mathfrak V_2)_\xi$ belong to the same 
one-dimensional space $\wedge^{n-l}V_\xi\subset \wedge^{n-l}\q^*$.

Thus,  $\mathfrak V_1$ and $\mathfrak V_2$
are two elements of a free $\bbk[\q^*]$-module (the module of regular sections of 
$\wedge^{n-l}T(\q^*)$), which is isomorphic to $\wedge^{n-l}\q^*\otimes \bbk[\q^*]$. 
Furthermore, 
$(\mathfrak V_1)_\xi$ and $(\mathfrak V_2)_\xi$ are linearly dependent 
as elements of $\wedge^{n-l}\q^*$ for any $\xi\in S$. 
It then follows that 
$\mathfrak V_1$ and $\mathfrak V_2$ are linearly dependent as elements of the
vector space $\wedge^{n-l}\q^*\otimes \bbk(\q^*)$ over the field $\bbk(\q^*)$.
Hence there are mutually prime  $F_1,F_2\in \bbk[\q^*]$ such that 
$F_1\mathfrak V_1=F_2\mathfrak V_2$. If $F_2$ were non-constant, then the section
$\mathfrak V_1$ would vanish on a divisor, which contradicts the codim--2 property.
Therefore, we may assume that $F_2\equiv 1$.

The equality $F_1\mathfrak V_1=\mathfrak V_2$ shows that 
$\deg\mathfrak V_1 \le \deg\mathfrak V_2$, that is,
$(n-l)/2 \le \sum_{i=1}^l\deg (f_i-1)$, which yields (i). 

If $\sum_{i=1}^l\deg f_i=(n+l)/2$, then $\deg F_1=0$, i.e., $F_1$ is 
a nonzero constant.
Therefore $\q^*_{reg}=S$. Since $\codim (\q^*\setminus S)\ge 2$, 
Theorem~\ref{skr} below and the fact that $\trdeg\bbk[\q^*]^Q=l$ 
guarantee us that $\bbk[\q^*]^Q=\bbk[f_1,\dots,f_l]$.
\end{proof}%
The following general result appears in \cite[Theorem\,1.1]{prya}. Its prototype is a theorem of 
Skryabin on algebras of invariants in a positive characteristic \cite[Theorem\,5.4]{skr}.
\begin{s}{Theorem}  \label{skr}
Let $V$ be a $\bbk$-vector space. Suppose that homogeneous polynomials 
$f_1,\dots,f_m\in\bbk[V]$
satisfy the property that $\codim_V \{v\in V\mid (\textsl{d}f_1)_v\wedge \ldots
\wedge (\textsl{d}f_m)_v=0\} \ge 2$. Then any $f\in\bbk[V]$ that is algebraic over
the subalgebra $\bbk[f_1,\dots,f_m]$ is necessarily contained in\/ $\bbk[f_1,\dots,f_m]$.
\end{s}%
\vskip-1ex
\begin{rem}{Remarks}   \label{coord-ner}
1. In the spirit of \cite{or}, Theorem~\ref{ner-vo} can be stated in the more general 
context of polynomial Poisson algebras and their centres.

2. The equality $F_1{\cdot}(\underbrace{\pi\wedge\pi\wedge\cdots\wedge \pi}_{{(n-l)/2}})=
(\textsl{d}f_1\wedge\ldots\wedge \textsl{d}f_l)^\diamond$ 
can be expressed in the coordinate form
as follows. Let $x_1,\dots,x_n$ be a basis for $\q$.
Form the $n\times n$ matrix 
$\Pi=([x_i,x_j])$ with entries in $\q=\bbk[\q^*]_1$. It is nothing but the matrix of 
the Poisson tensor $\pi$. If $I\subset [n]$ and $\#I=n-l$, then $\Pi_I$ 
denotes the Pfaffian of the principal $n{-}l$ submatrix of $\Pi$
corresponding to $I$. That is, $\Pi_I=\text{Pf}\bigl( ([x_i,x_j])_{i,j\in I}\bigr)$, and
it is a polynomial of degree $(n-l)/2$.
Another ingredient is the $l\times n$ matrix $\mathcal D=(\partial f_i/\partial x_j)$
of all partial derivatives of the polynomials $f_1,\dots, f_l$. Given $I$ as above, set
$\ov{I}=[n]\setminus I$. Let $\mathcal D_I$ denote the $l$ minor of $\mathcal D$
whose set of columns is $\ov{I}$. Then we have
\[
            F_1\Pi_I=\mathcal D_I \text{ for any } I\subset [n] 
              \text{ with } \#I=n-l \ . 
\]
Similar (although more complicated) 
equalities for minors were obtained in \cite[\S\,1]{rs1} for semisimple Lie algebras.

3. The proof of Theorem~\ref{ner-vo} shows that the equality
$F_1{\cdot}(\underbrace{\pi\wedge\pi\wedge\cdots\wedge \pi}_{{(n-l)/2}})=F_2{\cdot}
(\textsl{d}f_1\wedge\ldots\wedge \textsl{d}f_l)^\diamond$ with some $F_1,F_2
\in\bbk[\q^*]^Q$ holds for any Lie algebra. 
This allows to draw different conclusions under different assumptions. For instance,
if $\q$ is arbitrary and $f_1,\dots,f_l\in \bbk[\q^*]^Q$ have the property that $l=\ind\q$ 
and $S$ is big, then $\sum_i\deg f_i\le (\dim\q+\ind\q)/2$. 
\end{rem}%
{\bf Example}.  Let $\h=\bbk a+\bbk b +\bbk h$ be a Heisenberg Lie algebra
($[a,b]=h$ and $h$ is central). Here $\ind\h=1$ and $\bbk[\h^*]^H$ is generated
by $f_1=h$. Hence $1=\deg f_1< (\dim\h+\ind\h)/2=2$. This means that 
$\h$ does not have the codim--2 property, which is also easily verified directly.


\section{Semi-direct products, isotropy contractions, and $\BZ_2$-gradings}
\label{sect:recall}
\setcounter{equation}{0}

\noindent
Let $Q$ be a connected algebraic group with Lie algebra $\q$.

{\bf  (A)} {\it\un{Semi-direct} p\un{roducts}.} 
Let $V$ be a (finite-dimensional rational) $Q$-module, and hence a $\q$-module.
Then $\q\times V$ has a natural structure of Lie algebra,
$V$ being an Abelian ideal in it. Explicitly, if $x,x'\in \q$ and $v,v'\in V$, then
\[
   [(x,v), (x',v')]=([x,x'], x{\cdot}v'-x'{\cdot}v) \ .
\]
This Lie algebra is denoted by $\q\ltimes V$.
A connected algebraic group with Lie algebra $\q\ltimes V$ 
is identified set-theoretically with $Q\times V$, and 
we write $Q\ltimes V$ for it. The product in $Q\ltimes V$ is given by
\[
    (s,v)(s',v')= (ss', (s')^{-1}{\cdot}v+v') \ .
\]
In particular,  $(s,v)^{-1}=(s^{-1}, -s{\cdot}v)$.
The adjoint representation of $Q\ltimes V$ is given by the formula
\begin{equation}    \label{adj-QV}
(\Ad(s,v))(x',v')=(\Ad(s)x', s{\cdot}v'-x'{\cdot}v) \ ,
\end{equation}
where $v,v'\in V$, $x\in\q$, and $s\in Q$.
\\[.7ex]
Note that $V$ can be regarded as either a commutative
unipotent subgroup of $Q\ltimes V$ or a commutative nilpotent subalgebra of
$\q\ltimes V$. Referring to $V$ as subgroup of $Q\ltimes V$, 
we write $1\ltimes V$.

Set $\ka=\q\ltimes V$ and $K=Q\ltimes V$.
The dual space $\ka^*$ is identified with $\q^*\oplus V^*$, and
a typical element of it is denoted by $\eta=(\ap,\xi)$. 
The coadjoint representation of $\ka$ is given by
\begin{equation}  \label{action*}
   (\ads(x,v))(\ap,\xi)=(\ads(x)\ap - v\ast\xi, x{\cdot}\xi) \ . 
\end{equation}
Here the mapping $((x,\xi)\in \q\times V^*)\mapsto (x{\cdot}\xi\in V^*)$ is the 
natural $\q$-module structure on $V^*$, and 
$((v,\xi)\in V\times V^*)\mapsto (v\ast\xi\in \q^*)$ is the moment mapping
with respect to the symplectic structure on $V\times V^*$.

\begin{s}{Theorem}    \label{semidir_gen}
Let\/ $\ka=\q\ltimes V$ be an arbitrary semi-direct product. Then
\begin{itemize}
\item[\sf 1.] \ The algebras\/ $\bbk[\ka]^{K}$
and\/ $\bbk[\ka^*]^{K}$
are bi-graded;
\item[\sf 2.] \ There are natural inclusions 
$i_\q: \bbk[\q]^Q\hookrightarrow \bbk[\ka]^{K}$ and
$i_{V^*}:\bbk[V^*]^Q\hookrightarrow \bbk[\ka^*]^{K}$.
\item[\sf 3.] \ Let $J_1\subset \bbk[\ka]^{K}$
be the ideal of all
bi-homogeneous polynomials having a positive degree with respect to $V$.
Then $\bbk[\ka]^{K}=i_\q(\bbk[\q]^Q)\oplus J_1$;
\item[\sf 4.] \ Let $J_2\subset \bbk[\ka^*]^{K}$ be the ideal of all
bi-homogeneous polynomials having a positive degree with respect to $\q^*$.
Then $\bbk[\ka^*]^{K}=i_{V^*}(\bbk[V^*]^Q)\oplus J_2$.
\end{itemize}
\end{s}\begin{proof}
1. Let $\bbk[\ka^*]_{(a,b)}$ denote the space of bi-homogeneous polynomials of
degree $a$ (resp. $b$) with respect to $\q^*$ (resp. $V^*$).
Clearly, each $\bbk[\ka^*]_{(a,b)}$ is $Q$-stable. 
Given $v\in V$, let $D_{v,\ka^*}$ denote the derivation of $\bbk[\ka^*]$
corresponding to $(0,v)\in\ka$. Then Eq.~\eqref{action*} shows that 
$D_{v,\ka^*}(\bbk[\ka^*]_{(a,b)})\subset \bbk[\ka^*]_{(a-1,b+1)}$. Hence if 
$f\in \bbk[\ka^*]^{K}$ is a homogeneous polynomial, then all its bi-homogeneous
components are $Q$-invariant and the 
bi-homogeneous component of highest degree with respect to $V^*$ is  also 
$1\ltimes V$-invariant. Then we argue by induction.
\\[.6ex]
The similar argument works for $(\ka,\ad)$. Here 
$D_{v,\ka}(\bbk[\ka]_{(a,b)})\subset \bbk[\ka]_{(a+1,b-1)}$ and 
one has to consider the bi-homogeneous
component of $f\in\bbk[\ka]^{K}$ having the maximal degree
with respect to $\q$.

2. We can regard $\q$ and $V^*$ as $Q\ltimes V$-modules with trivial action of
$1\ltimes V$. Then
consider the natural surjective homomorphisms of $Q\ltimes V$-modules
$\q\ltimes V\to \q$ and $(\q\ltimes V)^* \to V^*$.

3,4. Obvious.
\end{proof}%
We will omit the indication of $i_\q$ and $i_{V^*}$ in the sequel. 
If $\q=\g$ is a reductive (algebraic) Lie algebra, then
$\textit{Inv}(\g\ltimes V,\ad)$
is always polynomial \cite[Theorem\,6.2]{p05}. This is, however, not always the case for
$\textit{Inv}(\g\ltimes V,\ads)$. For, it follows from 
Theorem~\ref{semidir_gen}(4) that any minimal generating system of $\bbk[V^*]^G$ is a part
of a minimal generating system of $\textit{Inv}(\g\ltimes V,\ads)$. In particular, if
$\textit{Inv}(\g\ltimes V,\ads)$ is polynomial, then so is $\bbk[V^*]^G$.

{\bf (B)}  {\it\un{Isotro}py \un{contractions}.}
Let $\h$ be a subalgebra of $\q$ such that $\q=\h\oplus\mathfrak m$ 
for some $\ad\h$-stable subspace $\mathfrak m\subset \q$.
(Such an $\h$ is said to be {\it reductive in\/} $\q$.) Then $\ma$ is an $H$-module.
If $\h$ is the fixed-point subalgebra of an involutory automorphism of $\q$, then
it is reductive in $\q$. In this case, $\h$ is called a {\it symmetric subalgebra} of $\q$
and the $(\q,\h)$ is called a {\it symmetric pair}.

\begin{rem}{Definition}   \label{def:degener}
If $\h$ is reductive in $\q$, then the representation of $H$ on $\ma$ is called the 
{\it isotropy representation\/} and the Lie algebra 
$\h\ltimes\mathfrak m$ is called an {\it isotropy contraction\/} 
of $\q$. If $\h$ is symmetric,
so the decomposition $\q=\h\oplus\mathfrak m$ is a $\mathbb Z_2$-grading,
then $\h\ltimes\mathfrak m$ is also called a $\BZ_2$-{\it contraction\/}
of $\q$.
\end{rem}%
Here $\h\ltimes\ma$ is a contraction of $\q$
in the  sense of the deformation theory of Lie algebras, see e.g.
\cite[Chapter 7, \S\,2]{t41}.  
More precisely, consider the invertible
linear map $c_t: \q\to \q$, $t\in \bbk\setminus\{0\}$, 
such that $c_t(h+m)=h+t^{-1}m$ \ ($h\in\h$, $m\in\mathfrak m$).  
Define the new Lie algebra multiplication $[\ ,\ ]_{(t)}$
on the vector space $\q$ by the rule
\[
     [x,y]_{(t)}:= c_t\bigl( [ c_t^{-1}(x), c_t^{-1}(y)]\bigr), \quad x,y\in\q \ .
\]
Then the algebras $\q_{(t)}$ are isomorphic for all $t\ne 0$, and
$\lim_{t\to 0}\q_{(t)}=\h\ltimes\mathfrak m$.

Suppose $\q=\g$ is reductive and $\h\subset\g$ is also reductive. Then
the isotropy representation of $H$ is orthogonalisable
(in particular, $\ma\simeq\ma^*$ as $H$-module) and 
$\ka:=\h\ltimes\ma$ is called a {\it reductive isotropy contraction\/} (of $\g$). Here
$K^u:=1\ltimes \ma$ is the unipotent radical of $K=H\ltimes \ma$.

A natural hope is that the algebras $\bbk[\ka]^K$ and $\bbk[\ka^*]^K$
could keep some good properties of $\bbk[\g]^G$. But this is not always the case.
For instance, the transcendence degree of 
$\bbk[\ka]^K$ and $\bbk[\ka^*]^K$ can be larger than $\rk\g$, and to guarantee equalities,
one has to impose different constraints on $\h$.
Since $\h$ is reductive,  $\bbk[\ka]^K$ is polynomial \cite[Theorem\,6.2]{p05}; 
in other words, $\ka\md K$ is an affine space.
For future reference, we record the following fact.
\begin{s}{Proposition {\ququ \cite[Prop.\,9.3]{p05}}}  \label{9.3} 
Let $\ka$ be a reductive isotropy contraction of $\g$. Then 

1. \ $\dim\ka\md K=\rk\g$ if and only if\/ $\h$ contains a regular semisimple element of $\g$;

2. \ $\ind\ka=\rk\g$ if and only if\/ $G/H$ is a spherical homogeneous space.
\end{s}%
Both these conditions are satisfied if $\h$ is a symmetric subalgebra of $\g$.
However, the adjoint and coadjoint representations of $\ka$ are quite different, and
should be studied separately.

\begin{s}{Lemma}   \label{stab}
If\/ $\ka=\h\ltimes\ma$ is a reductive isotropy contraction of 
$\g$, then the quotient field of\/ $\bbk[\ka^*]^K$ equals\/ $\bbk(\ka^*)^K$.
\end{s}\begin{proof}
We have $K=T_H{\cdot}(K,K)$, where $T_H$ is connected centre of $H$ and $(K,K)$ is
the derived group of $K$. Since $(K,K)$ has no rational character, 
the quotient field of $\bbk[\ka^*]^{(K,K)}$ equals $\bbk(\ka^*)^{(K,K)}$.
It follows that any $f\in \bbk(\ka^*)^K$ can be written as $f=f_1/f_2$, where 
$f_1,f_2\in  \bbk[\ka^*]^{(K,K)}$ are semi-invariants of $T_H$ of the same weight,
say $\chi$. Clearly, if $\bbk[\ka^*]$ contains a semi-invariant of $T_H$ of weight
$\nu$, then it also contains a semi-invariant of weight $-\nu$. (Because $\ka^*\simeq\g$
as $T_H$-modules.) The same assertion is also true for $\bbk[\ka^*]^{(K,K)}$ in  place
of $\bbk[\ka^*]$. [Use the fact that the automorphism of $K$ (as a variety!) that is trivial on $(K,K)$ and
takes $t$ to $t^{-1}$ for any $t\in T_H$ does not change the $K$-action on $\ka^*$.]
Thus, if $h\in  \bbk[\ka^*]^{(K,K)}$ is a semi-invariant of weight $-\chi$, then
$f=(f_1h)/(f_2h)$, and we are done.
\end{proof}%
{\bf (C)} {\it $\BZ_2$-g\un{radin}g\un{s} 
\un{o}f \un{reductive} \un{Lie} \un{al}g\un{ebras}.}
Let $G$ be a connected reductive algebraic group with $\g=\text{Lie\,}G$. 
Let $\N$ denote the set of nilpotent elements of $\g$. If $\g=\g_0\oplus\g_1$ is a 
$\BZ_2$-grading of $\g$, then $G_0$ is the connected subgroup of 
$G$ with Lie algebra $\g_0$. Recall some results on 
the isotropy representation $(G_0:\g_1)$. The standard reference for this is \cite{kr}.

-- Any $v\in\g_1$ admits a unique decomposition $v=v_s+v_n$, where $v_s\in\g_1$ is 
semisimple and $v_n\in \N\cap\g_1$; $v=v_s$ if and only if $G_0{\cdot}v$ is closed;
$v=v_n$ if and only if the closure of $G_0{\cdot}v$ contains the origin.
For any $v\in\g_1$, there is the induced $\BZ_2$-grading of the centraliser 
\ $\g_v=\g_{0,v}\oplus \g_{1,v}$, and $\dim\g_0-\dim\g_{0,v}=\dim\g_1-\dim\g_{1,v}$.

-- Let $\ce\subset\g_1$  be a maximal subspace consisting of pairwise commuting
semisimple elements. Any such subspace is called a {\it Cartan
subspace\/}. All Cartan subspaces are $G_0$-conjugate and $G_0{\cdot}\ce$ is dense in $\g_1$; 
$\dim\ce$ is called the {\it rank\/} of the $\BZ_2$-grading or pair $(\g,\g_0)$, denoted
$\rk(\g,\g_0)$.  If $v\in\ce$ is $G_0$-regular (i.e.,
$\dim G_0{\cdot}v$ is maximal), then $\g_{1,v}=\ce$ and $\g_{0,v}$ is a generic stabiliser
for the $G_0$-module $\g_1$.

-- The algebra $\bbk[\g_1]^{G_0}$ is polynomial and $\dim\g_1\md G_0=\rk(\g,\g_0)$. 
The quotient map $\pi: \g_1\to \g_1\md G_0$ is equidimensional.
We write $\gN(\g_1)$ for $\pi^{-1}(\pi(0))$.
Any fibre of $\pi$ contains finitely many $G_0$-orbits and each closed $G_0$-orbit in
$\g_1$ meets $\ce$. There is a finite reflection group $W_\ce\subset GL(\ce)$ 
("the little Weyl group") such that $\ce/W_\ce \simeq \g_1\md G_0$.

{\bf (D)} {\it \un{Reductive} $\BZ_2$-\un{contractions}.}
Given a $\BZ_2$-grading of $\g$,
consider the $\BZ_2$-contraction $\ka=\g_0\ltimes\g_1$.
Set $K^u:=1\ltimes \g_1$. 
The adjoint representation of $\ka$ was studied in \cite{p05}.
Below we summarise the relevant invariant-theoretic results, see \cite
[Prop.\,5.3 \& Theorem\,9.13]{p05}:

\begin{itemize}
\item \  Let $\te_0$ be a Cartan subalgebra of $\g_0$ and $\g_1^0$ the centraliser of $\te_0$
in $\g_1$. Then $\te_0\ltimes \g_1^0$ is a generic stabiliser for $(\ka,\ad)$.
[As is well-known, $\g_0$ contains regular
semisimple elements of $\g$. Therefore $\dim(\te_0\ltimes\g_1^0)=\rk\g$.]
\item \  $\bbk[\ka]^{K^u}$ is a polynomial algebra of Krull dimension $\dim\g_0+\dim\g_1^0$;
\item \  $\bbk[\ka]^K$ is a polynomial algebra of Krull dimension $\rk\g$;
\item \  the quotient map $\pi_\ka:\ka\to \ka\md K$ is equidimensional and $\bbk[\ka]$ is a free $\bbk[\ka]^K$-module;
\item \  ${\goth N}(\ka):=\pi_\ka^{-1}(\pi_\ka((0))$ is an irreducible complete intersection. 
       If\/ $\bbk[\ka]^K=\bbk[f_1,\ldots, f_l]$, $l=\rk\g$, then 
       the ideal of ${\goth N}(\ka)$ in $\bbk[\ka]$ is generated by $f_1,\ldots, f_l$.
\end{itemize} 
However, the key fact is that 
there is a natural description of basic invariants in $\bbk[\ka]^K$ 
(see \cite[Sect.\,6]{p05}), which enables us to prove the above results.
Namely, the set of basic invariants consists of two parts.
First, we take a set of basic invariants in $\bbk[\g_0]^{G_0}$,
say $f_1,\dots,f_m$. Here $m=\rk\g_0$. Next, we consider the
set, $\Mor_{G_0}(\g_0,\g_1)$, of all $G_0$-equivariant polynomial
morphisms $\tau:\g_0\to\g_1$. By a result of Kostant \cite{ko63},
$\Mor_{G_0}(\g_0,\g_1)$ is a free $\bbk[\g_0]^{G_0}$-module of rank $\dim\g_1^0=l-m$.
Given $F\in \Mor_{G_0}(\g_0,\g_1)$, define the polynomial 
$\widehat F\in \bbk[\ka]$ by $\widehat F(x_0,x_1)=\langle F(x_0),x_1\rangle$.
Here $x_i\in\g_i$ and $\langle\ ,\ \rangle$ stands for a nondegenerate $G_0$-invariant
symmetric bilinear form on $\g_1$. It is easily seen that $\widehat F\in \bbk[\ka]^K$.
If $F_1,\dots,F_{l-m}$ is a basis for $\Mor_{G_0}(\g_0,\g_1)$,
then $f_1,\dots,f_m,\widehat F_1,\dots,\widehat F_{l-m}$ is a set of basic invariants
in $\bbk[\ka]^K$.

{\bf Remark}.
It seems that the reason for success in case of $(\ka,\ad)$ 
is that $\g_0$ always contains regular semisimple 
elements of $\g$. We will see in Section~\ref{sect:N-reg} that if $\g_1$ contains a 
regular semisimple element, then
$\bbk[\ka^*]^K$ is polynomial and, moreover, there is a similar description of basic invariants
and similar properties hold.


\section{Constructing invariants for reductive $\BZ_2$-contractions}
\label{sect:z2-contr}
\setcounter{equation}{0}

\noindent
From now on, $\ka=\g_0\ltimes\g_1$ is a reductive $\BZ_2$-contraction and $K=G_0\ltimes \g_1$. 
Our primary goal is to study invariant-theoretic properties of
the coadjoint representation of $\ka$. However, we also mention 
results for $(\ka,\ad)$, if they are parallel to those for $(\ka,\ads)$ and
are not contained in \cite{p05}.

We identify the $G$-modules $\g$ and $\g^*$, using a non-degenerate
$G$-invariant symmetric bilinear form on $\g$.
Moreover, $\g_i$ and  $\g^*_i$ (and hence $\ka$ and $\ka^*$) are identified 
as $G_0$-modules. This means, for instance, that we can speak about a Cartan subspace of
$\g_1^*$ and that any $f\in\bbk[\g]$ can also be regarded as function on
$\ka$ or $\ka^*$. Usually, it is clear from the context whether $\g_i$ 
is regarded as a subspace of $\g$ or $\ka$ or $\ka^*$. (This makes no difference
as long as only $G_0$-module structure is involved.) However, if we wish to stress
that $\g_i$ is regarded as subspace of $\ka^*$, then we  write $\g_i^*$ for it.

There is a natural procedure of getting elements of $\bbk[\ka]^K$ and $\bbk[\ka^*]^K$
via "$\BZ_2$-dege\-ne\-ra\-tions" of $G$-invariants on $\g$.
Let $\bbk[\g]_{(a,b)}$ denote the space of bi-homogeneous polynomials
of degree $a$ with respect to $\g_0$ and degree $b$ with respect to $\g_1$.

Given a homogeneous polynomial $f\in \bbk[\g]$ of degree $n$, 
let us decompose $f$ into the sum of bi-homogeneous components 
$f=\sum_{i=k}^m f_i$, where $f_i\in \bbk[\g]_{(n-i,i)}$ and
it is assumed that $f_k, f_m\ne 0$. Then we set
$f^\bullet:=f_m$ and $f_\bullet:=f_k$.

\begin{s}{Proposition}  \label{z2-contra}
Suppose that $f\in \bbk[\g]^G$ is homogeneous. Then
\begin{itemize}
\item[\sf (i)] \ regarding $f$ as function on $\ka$, we have 
$f_\bullet\in \bbk[\ka]^K$;
\item[\sf (ii)] \ regarding $f$ as function on $\ka^*$, we have
$f^\bullet\in \bbk[\ka^*]^K$.
\end{itemize} 
\end{s}\begin{proof}
Clearly, each $f_i$ is $G_0$-invariant. The derivation of $\bbk[\g]$ corresponding
to $x\in\g_1$ is denoted by $D_{x,\g}$. The commutator relations for the $\BZ_2$-grading
show that for $f_i\in \bbk[\g]_{(n-i,i)}$, we have 
\[
D_{x,\g}(f_i)\in \bbk[\g]_{(n-i-1,i+1)}\oplus \bbk[\g]_{(n-i+1,i-1)} \ .
\]
Accordingly, we write $D_{x,\g}=D_x^{(+1)}+D_x^{(-1)}$, where
$D_x^{(+1)}: \bbk[\g]_{(n-i,i)} \to \bbk[\g]_{(n-i-1,i+1)}$. It follows
that $D_x^{(+1)}(f_m)=0$ and $D_x^{(-1)}(f_k)=0$, if $f\in\bbk[\g]^G$.
The key observation is that 
$D_x^{(+1)}=D_{x,\ka^*}$ and $D_x^{(-1)}=D_{x,\ka}$.
\end{proof}%
Part (ii) appears in \cite{bra}. However, for some particular cases, this construction of 
invariants of the coadjoint representation is considered in \cite{rosen}.
The passages $f\mapsto f^\bullet$ and $f\mapsto f_\bullet$ will be referred to as
$\BZ_2$-{\it degenerations\/} of (homogeneous) invariants in $\bbk[\g]^G$.
In this way, one obtains bi-graded subalgebras 
\[
  \textrm{gr}_\bullet(\bbk[\g]^G):=\{f_\bullet\mid f\in\bbk[\g]^G\}\subset \bbk[\ka]^K      
      \text{ and }
  \textrm{gr}^\bullet(\bbk[\g]^G):=\{f^\bullet\mid f\in\bbk[\g]^G\}\subset \bbk[\ka^*]^K \ .
\]
However, both inclusions can be strict. For, the $\BZ_2$-degeneration preserves 
the usual degree of polynomials, but it is possible in many cases
to point out an element of $\bbk[\ka]^K$ or $\bbk[\ka^*]^K$ whose degree
does not occur as degree of elements of $\bbk[\g]^G$.
For instance, if $\rk\g_0=\rk\g$, then $\bbk[\ka]^K\simeq \bbk[\g_0]^{G_0}$.
Clearly, $\bbk[\g_0]^{G_0}$ has "more" elements than $\bbk[\g]^G$.
Examples for $\bbk[\ka^*]^K$ are discussed in Remark~\ref{bad}.

\begin{rem}{Remark}   \label{covar*}
As is explained in Section~\ref{sect:recall}(D),
invariants of $(\ka,\ad)$ can be constructed using the $\bbk[\g_0]^{G_0}$-module
({\it module of covariants}) $\Mor_{G_0}(\g_0,\g_1)$. One might suggest that there was a similar procedure for 
$(\ka,\ads)$, which makes use of the module of covariants $\Mor_{G_0}(\g_1^*,\g_0)$.
However, this does not always work.
For $F\in \Mor_{G_0}(\g_1^*,\g_0)$, we can define $\widehat F\in \bbk[\ka^*]$ by 
 \ $\widehat F(\xi_0,\xi_1)=\langle F(\xi_1),\xi_0\rangle$,
where $\xi_i\in\g_i^*$ and $\langle\ ,\ \rangle$ is a $G_0$-invariant
non-degenerate symmetric bilinear form on $\g_0$.
Obviously, $\widehat F$ is $G_0$-invariant. But its invariance relative to 
$K^u=1\ltimes\g_1$ reduces to the condition that
\[
       F(\xi) \in \g_{0,\xi}  \text{ \ for all } \xi\in\g_1^* \ .
\]
This condition and $G_0$-equivariance of $F$ show that
$F(\xi)$ belong to the centre of $\g_{0,\xi}$. That is, such a nonzero covariant
may only exist if a generic stabiliser for the $G_0$-module $\g_1^*$ has a non-trivial
centre. 
\end{rem}%
\vskip-1.5ex
\begin{s}{Theorem}   \label{thm:codim2}
Any reductive $\BZ_2$-contraction  has the codim-2 property for\/ $\ads$.
\end{s}\begin{proof}
(a) We explicitly describe certain big open subset of $\ka^*$ that is
contained in $\ka^*_{reg}$. 

Let $\eta=(\ap,\xi)\in \ka^*$ be an arbitrary point, 
where $\ap\in\g_0^*$ and $\xi\in\g_1^*$. 
Write $\g_{0,\xi}$ for the stabiliser of $\xi$ in $\g_0$.
Then $\g_1\ast\xi=\Ann(\g_{0,\xi})\subset\g_0^*$ and 
therefore $\g_0^*/\g_1\ast\xi \simeq \g_{0,\xi}^*$.
Using the last isomorphism, we let $\bar\ap$ denote the image of $\ap$ in $\g_{0,\xi}^*$.
By  \cite[Prop.\,5.5]{p05}, 
\begin{equation}  \label{eq:dim_stab}
\dim \ka_{\eta}=\codim_{\g_1^*}(G_0{\cdot}\xi)+\dim (\g_{0,\xi})_{\bar\ap} \ ,
\end{equation}
where the last summand refers to the stabliser of $\bar\ap$ with respect to the
coadjoint representation of $\g_{0,\xi}$.

Let $\Omega\subset \g_1^*$ be the open subset of $G_0$-regular points, i.e.,
\[
\Omega=\{\xi\in\g_1^*\mid \dim G_0{\cdot}{\xi}=\dim\g_1-\rk(\g,\g_0)\}
\ .
\] 
It follows from Eq.~\eqref{eq:dim_stab} that in order to obtain a $K$-regular point in $\ka^*$,
it suffices to take a $G_0$-regular point $\xi\in\g_1^*$ and then, if the
equality $\ind\g_{0,\xi}=\rk\g-\rk(\g,\g_0)$ holds, 
to take an $\ap$ such that 
$\bar\ap\in\g_{0,\xi}^*$ is a $G_{0,\xi}$-regular point.
Let us prove that the set of such points $(\ap,\xi)$ contains a big open subset of $\ka^*$.

Let $\pi:\g_1^*\to \g_1^*\md G_0$ be the quotient mapping.
Consider the Luna stratification of $\g_1^*\md G_0$ \cite[III.2]{mem33}.
(Recall that $\nu,\nu'\in \g_1^*\md G_0$ belong to the same stratum, if 
the closed $G_0$-orbits in $\pi^{-1}(\nu)$ and $\pi^{-1}(\nu')$
are isomorphic as $G_0$-varieties.) An exposition of Luna's theory can also be found in
\cite{slo}.
Let $(\g_1^*\md G_0)_i$ be the union of all strata of codimension $i$.
For instance, $(\g_1^*\md G_0)_0$ is the unique open stratum.
Set $\Omega_i=\pi^{-1}((\g_1^*\md G_0)_i)\cap \Omega$. 
Since $\pi$ is equidimensional and each fibre of $\pi$ meets $\Omega$, 
$\codim_{\g_1^*}\Omega_i=i$.
In particular, $\Omega_0\cup\Omega_1$ is a big open subset of $\g_1^*$ and hence
$(\Omega_0\cup\Omega_1)\times \g_0^*$ is a big open subset of $\ka^*$.
Let us prove that $\ka^*_{reg}\cap((\Omega_0\cup\Omega_1)\times \g_0^*)$ is still big.

If $\xi\in\Omega_0$, then $\xi$ is semisimple and $\g_{0,\xi}$ is reductive. 
Since $(\g_{0,\xi},\ads)$ has codim--2 property, the set  
$\ka^*_{reg}\cap (\Omega_0\times \g_0^*)$ is big in $\Omega_0\times \g_0^*$
(but not in $\ka^*$ !). To obtain a big subset of $\ka^*$, we have to check that
$\ka^*_{reg}\cap (\Omega_1\times \g_0^*)$ is dense in $\Omega_1\times \g_0^*$.
In view of the previous discussion, this amounts to the verification of the 
equality $\ind\g_{0,\xi}=\rk\g-\rk(\g,\g_0)$ for any $\xi\in \Omega_1$.

Using the Jordan decomposition in $\g_1^*$ and taking the centraliser of the
semisimple part of $\xi\in\Omega_1$, one reduces the problem to the case of 
symmetric pairs of rank 1. Namely, let $\xi=\xi_s+\xi_n$, where the semisimple 
element $\xi_s$ belong to a fixed Cartan subspace $\ce$.
Then the centraliser of
$\xi_s$ in $\g$ has the following structure: $\z_\g(\xi_s)=\ah\dotplus\h$,
where $\ah\subset\ce$, $\dim\ah=\dim\ce-1$, $\h$ is reductive, and the induced
$\BZ_2$-grading of $\h$ has rank 1. Furthermore, $\xi_n\in\h_1\subset\h$ and
$\g_{0,\xi}=\h_{0,\xi_n}$.  Hence it remains to handle the rank one case.

(b) \ Suppose $\rk(\g,\g_0)=1$, i.e., $\dim\g_1^*\md G_0=1$. 
Then
$(\g_1^*\md G_0)_1=\{pt\}=\pi(0)$ and $\Omega_1$ is the set of 
$G_0$-regular nilpotent elements of $\g_1$.
Here we have to check that if $\xi\in\Omega_1$,
then $\ind\g_{0,\xi}=\rk\g-1$. Since $\g_\xi=\bbk\xi\dotplus \g_{0,\xi}$ (a direct sum
of Lie algebras), we need actually the equality $\ind\g_\xi=\rk\g$. Such an equality is known
as "Elashvili's conjecture", and it is proved for all $\xi$ in the classical Lie algebras in 
\cite{ksana}.
The only non-classical symmetric pair of rank one is $(\GR{F}{4},\GR{B}{4})$,
where one has to test the stabiliser of a sole nilpotent $G_0$-orbit.  
Here the isotropy representation is the spinor representation
of $\GR{B}{4}$. By Igusa's computations \cite{igusa}, the stabiliser 
$\g_{0,\xi}$ is the semi-direct product of $\GR{G}{2}$ and its 7-dimensional representation.
Then using Ra\"\i s' formula \cite{rais}, we obtain $\ind\g_{0,\xi}=3$.
(It is also easy to perform similar verifications for the three classical
series of symmetric pairs of rank one.)
\end{proof}%
Combining Proposition~\ref{9.3}(2), Lemma~\ref{stab},
Theorems~\ref{ner-vo} and \ref{thm:codim2}, we obtain 

\begin{s}{Corollary}  \label{z2-sdeg}
If $f_1,\dots,f_l\in\bbk[\ka^*]^K$ are homogeneous algebraically independent 
and $l=\rk\g$, then $\sum_{i=1}^l \deg f_i \ge (\dim\g+l)/2$.
\end{s}%
As $\bbk[\ka^*]^K$ is a bi-graded algebra (Theorem~\ref{semidir_gen}), 
one can take bi-homogeneous polynomials $f_1,\dots,f_l$. 
Our next goal is to provide a "bi-graded" refinement of Corollary~\ref{z2-sdeg}.

For $f\in \bbk[\ka^*]_{(a,b)}$, we write $\bideg f=(a,b)$. 
Here $a$ and $b$ refer to  the $\g_0^*$-degree and $\g_1^*$-degree, respectively.
Let $\es\subset\g_0$ be a generic stabiliser for the isotropy representation $(G_0:\g_1)$.
It is a reductive Lie algebra. 

\begin{s}{Theorem}   \label{bideg}
Suppose that $f_1,\dots,f_l\in\bbk[\ka^*]^K$ are bi-homogeneous algebraically independent 
and $l=\rk\g$. Then 
$\sum_{i=1}^l\bideg f_i \ge ((\dim\es+\rk\es)/2, \dim\g_1)$ (componentwise).
\end{s}\begin{proof}
First of all, the inequality in question is a refinement of that in 
Corollary~\ref{z2-sdeg}. Indeed,  
\[
\dim\g_1\md G_0=\dim\g_1-\max_{x\in\g_1}\dim G_0{\cdot}x=\dim\g_1-\dim\g_0+\dim\es\ .
\]
On the other hand, $\dim\g_1\md G_0=\dim\ce=\rk\g-\rk\es$.
Equating two expressions for $\dim\g_1\md G_0$ and
rearranging them, we obtain  $\dim\g_1 +(\dim\es+\rk\es)/2 = (\dim\g+\rk\g)/2$, 
as required.

To prove the inequality, we use the construction of Theorem~\ref{ner-vo}
in the coordinate form, as described in Remark~\ref{coord-ner}(2).
Let us match $\es$ and $\ce$ such that $\es$ is the stabliser in $\g_0$ of a 
generic element  $x\in\ce$.
Then $\es\oplus\ce=\g_x$ is a Levi subalgebra.
Let \ $\te(\es)$ be a Cartan subalgebra of $\es$. By  \cite[Sect.\,5]{p05},
$\h=\te(\es)\ltimes\ce\subset\ka$ is a generic stabiliser for $(\ka,\ads)$.
We may (and will) consider $\h$ as a subspace in either $\ka$ or $\ka^*$.
In the last case we will denote it as $\h^*=\te(\es)^*\oplus\ce^*$.
In our situation, $\h$ has the property that $\z_\ka(\h)=\h$. It then follows from
\cite[Theorem\,3.4]{p05} that $(\ka^*)^\h=\h^*$ and $\ka=[\ka,\h]\oplus\h$. 
Taking the annihilators, we obtain the dual decomposition
 $\ka^*=\ads(\ka){\cdot}\xi\oplus \h^*$, where $\xi\in\h^*$ is a
generic point. 


Choose a basis $(x_1,\dots,x_n)$ for $\ka$ such that 
$(x_1,\dots,x_{n-l})$ is a basis for $[\ka,\h]$ and $(x_{n-l+1},\dots,x_n)$ is a basis
for $\h$. 
Recall that if $I\subset [n]$ and $\# I=n-l$, then 
$\Pi_I=\text{Pf}\bigl( ([x_i,x_j])_{i,j\in I}\bigr)$ and 
$D_I=\det \bigl((\partial f_i/\partial x_j)_{j\not\in I}\bigr)$.
Set $I_0=[n{-}l]$ and consider $\Pi_{I_0}$ and $\mathcal D_{I_0}$.
More precisely, we need the restriction of these polynomials to the subspace $\h^*$,
$\bar\Pi_{I_0}=\Pi_{I_0}\vert_{\h^*}$ and $\bar{\mathcal D}_{I_0}=
\mathcal D_{I_0}\vert_{\h^*}$.
Clearly, $\bar{\mathcal D}_{I_0}$ is the Jacobian of $f_1\vert_{\h^*},\dots, f_l\vert_{\h^*}$.
Hence $\bideg\bar{\mathcal D}_{I_0}=\bigl(\sum_{i=1}^l\bideg f_i\bigr)- (\rk\es,\dim\ce)$.

{\bf Claim~1.} $\bar\Pi_{I_0}\ne 0$.

{\sl Proof of Claim~1.} It is easily seen that $([x_i,x_j]\vert_{\h^*})_{i,j\in I}$
has a zero column unless $I=I_0$, hence $\Pi_I\vert_{\h^*}=0$ unless $I=I_0$.
The definition of  generic stabilisers says that
$K{\cdot}\h^*$ is dense in $\ka^*$. Since $\pi^{(n-l)/2}$ is $K$-invariant and 
the functions $\Pi_I$, $I\subset [n]$, are the coefficients of $\pi^{(n-l)/2}$ in
the basis $\{\wedge_{i\in I} x_i\mid I\subset [n]\}$,
they all cannot vanish on
$\h^*$.  \hfill $\diamond$

{\bf Claim~2.} $\bideg\bar\Pi_{I_0}=\bigl((\dim\es-\rk\es)/2,\dim\g_1-\dim\ce\bigr)$.

{\sl Proof of Claim~2.} Since $[\ka,\h]=(\h^*)^\perp$, this space is a sum of its 
intersections with $\g_0^*$ and $\g_1^*$. More precisely, using nondegenerate $G_0$-invariant 
symmetric bilinear forms on $\g_0$ and $\g_1$, we obtain
$[\ka,\h]_0=\te(\es)^\perp$ and $[\ka,\h]_1=\ce^\perp$. Hence 
$\dim [\ka,\h]_0=\dim\g_0-\rk\es$ and $\dim[\ka,\h]_1=\dim\g_1-\dim\ce$. 
Since $\dim\g_0-\dim\es=\dim\g_1-\dim\ce$, we have
$\dim [\ka,\h]_0=\dim[\ka,\h]_1 +(\dim\es-\rk\es)$. Assume that a basis 
for $[\ka,\h]$ is chosen such that we first have a basis for $\te(\es)^\perp\cap\es$,
then a basis for $\es^\perp\cap [\ka,\h]_0$, and finally a basis for $[\ka,\h]_1$. It is easily
seen that, for this choice of a basis, the matrix $([x_i,x_j]\vert_{\h^*})_{i,j\in I_0}$ is of
the form:
$\begin{pmatrix} A & 0    & 0 \\
                0 & \ast & B \\
                0 & -B^t & 0 \end{pmatrix}$, where $A$ is a skew-symmetric matrix of
order $\dim\es-\rk\es$, with entries in $\g_0$; $B$ is a square matrix of order 
$\dim\g_1-\dim\ce$, with entries in $\g_1$. It follows that 
$\bar\Pi_{I_0}=\text{Pf}(A)\det(B)$ has the required bi-degree. \hfill $\diamond$

By Remark~\ref{coord-ner}(2), there is an $F\in \bbk[\ka^*]$ such that 
$F\Pi_I=\mathcal D_I$ for any $I\subset [n]$. Applying this to $I_0$ shows
that $\bideg\bar\Pi_{I_0} \le \bideg \bar{\mathcal D}_{I_0}$, 
which completes the proof of theorem.
\end{proof}%
%


\section{Good generating systems for invariants associated with symmetric pairs}
\label{good_systems}
\setcounter{equation}{0}

\noindent
The presence of codim--2 property for the $\BZ_2$-contractions and the procedure of 
$\BZ_2$-degeneration of invariants enable us to state
a helpful sufficient condition for the polynomiality of $\bbk[\ka^*]^K$.

\begin{rem}{Definition}  \label{def:good}
Let $f_1,\dots,f_l\in \bbk[\g]^G$ be a set of basic invariants. 
We say that it is a {\it good generating system\/} for $(\g,\g_0)$ if the
$\BZ_2$-degenerations
$f_1^\bullet,\dots,f_l^\bullet\in \bbk[\ka^*]^K$ are algebraically independent.
\end{rem}%
\vskip-1ex

\begin{s}{Theorem}  \label{thm:good}
If $f_1,\dots,f_l$ is a good generating system for $(\g,\g_0)$,
then {\sf (i)} \  $\bbk[\ka^*]^K$ is freely generated by $f_1^\bullet,\dots,f_l^\bullet$ and
{\sf (ii)} \ $(\textsl{d}f_1^\bullet)_\xi,\dots,(\textsl{d}f_l^\bullet)_\xi$ are linearly independent
if and only if $\xi\in \ka^*_{reg}$. Furthermore, in this case, $\bbk[\ka^*]^K=
\gr^\bullet (\bbk[\g]^G)$.
\end{s}\begin{proof}
Since $\deg f_i=\deg f_i^\bullet$, $\rk\g=\ind\ka$, and $(\ka,\ads)$ has codim--2
property,  Theorem~\ref{ner-vo}(ii) applies to $f_1^\bullet,\dots,f_l^\bullet$.
\end{proof}%
The property of being `good' for a generating system is rather specific and can easily
be disturbed. For instance, if $\g=\sono$, then the coefficients of the characteristic 
polynomial of a matrix $M\in\sono$ form a good generating system for any symmetric
pair $(\sono,\mathfrak{so}_m\dotplus\mathfrak{so}_{2n+1-m})$ 
(see Theorem~\ref{good:so} below). 
But the polynomials $\tr(M^{2i})$, $i=1,\dots,n$, do not
form a good a generating system. 

\begin{rem}{Remark}  \label{bad}
Good generating systems do not always 
exist. For instance, consider the symmetric pair $(\GR{E}{6},\GR{F}{4})$.
Here $\bbk[\g_1]^{G_0}$ is freely generated by two polynomials of degree $2$
and $3$. Since $\bbk[\g_1]^{G_0}\hookrightarrow \bbk[\ka^*]^K$, the latter
has an element of degree three. However, the basic degrees of $\GR{E}{6}$ are
$2,5,6,8,9,12$. Hence $ \bbk[\g]^G$ does not contain elements of degree $3$
and the equality $\bbk[\ka^*]^K=\gr^\bullet (\bbk[\g]^G)$ cannot hold.
Similar phenomenon occurs for three other symmetric pairs:
$(\GR{E}{6},\GR{D}{5}\dotplus \te_1)$, $(\GR{E}{7},\GR{E}{6}\dotplus \te_1)$,
$(\GR{E}{8},\GR{E}{7}\dotplus\GR{A}{1})$. These are precisely the symmetric pairs
such that the restriction homomorphism $\bbk[\g]^G \to \bbk[\g_1]^{G_0}$ is not onto
\cite{helga}.
\end{rem}%
For some symmetric pairs,  it is possible 
to check directly that certain generating system is good. 
Below, we consider several examples.

{\sl \un{Practical} \un{tricks}.} 1. \ 
To prove that some polynomials in $\bbk[\ka^*]^K$ are algebraically independent,
it suffices to verify this for their restriction to a subspace. In case of $\BZ_2$-contractions,
it is convenient to take the subspace $\ce\oplus\es$, where $\ce$ is a fixed
Cartan subspace of $\g_1$ and $\es=\z_{\g_0}(\ce)$.  Recall that $\es$ is a generic stabiliser
for the $G_0$-module $\g_1$ and  $\es\oplus\ce=\g_x$ for a generic $x\in\ce$. 
Following our convention, we also regard $\es\oplus\ce$ as a subspace of $\ka^*$.  
Furthermore, one can work
with the smaller subspace $\ce\oplus\te(\es)$. 
Notice that this vector space has three masks:
as subspace of $\g$ it is just a Cartan subalgebra;
as a subspace of $\ka$ it is a generic stabiliser for $(\ka,\ads)$, say $\h$;
as a subspace of $\ka^*$ it is the fixed point space of $\h$.

2. \ Another useful observation is that if $f\in\bbk[\g]^G$, then taking the
restriction of $f^\bullet$ to $\ce\oplus\es$ (or $\ce\oplus\te(\es)$) 
is the same as first restricting  $f$ to 
$\ce\oplus\es$ (or $\ce\oplus\te(\es)$) 
and then taking the component of highest degree with respect to $\ce$. The reason is that
$f^\bullet\vert_{\ce\oplus\te(\es)}\ne 0$, since $f^\bullet\in \bbk[\ka^*]^K$ and
$\Ads K{\cdot}(\ce\oplus\te(\es))$ is dense in $\ka^*$.

\begin{s}{Theorem}  \label{good:so}
There is a good generating system for $(\g,\g_0)=
(\mathfrak{so}_{n+m},\mathfrak{so}_m\dotplus\son)$.
\end{s}\begin{proof}
Here $l=\lfloor (n+m)/2 \rfloor$ and $\rk(\g,\g_0)=\min\{n,m\}$.
We use the natural matrix model for $(\g,\g_0)$:

$\g_0=\left\{\begin{pmatrix} A & 0 \\ 0 & B \end{pmatrix}\right\}$ \ and \ 
$\g_1=\left\{\begin{pmatrix} 0 & C \\ -C^t & 0 \end{pmatrix}\right\}$,
\\[.6ex]
where $A$ (resp. $B$) is a skew-symmetric matrix of order $m$ (resp. $n$)
and $C$ is an $m\times n$ matrix. 
Assume that $m\le n$. Then $\es\simeq\mathfrak{so}_{n-m}$ and $\es\subset\son$.

For a Cartan subspace $\ce$, we take the set of matrices $C$ with only nonzero entries 
along the diagonal
starting in the upper left corner of $C$. Then $\es$ is the lower-right submatrix of $B$
of order $n-m$. That is, taking the partition of $M$ into the nine submatrices
corresponding to the sizes  $m,m,n-m$, we obtain 
$\ce\oplus\es=\left\{\tilde M=\begin{pmatrix} 0 & D & 0 \\  -D & 0 & 0 \\ 0 & 0 & E 
\end{pmatrix}\right\}$, where $D$ is a diagonal matrix of order $m$ and
$E$ is a skew-symmetric matrix of order $n-m$. Let $d_1,\dots,d_m$ be the diagonal
entries of $D$.

For a skew-symmetric matrix  $M$, let $f_i(M)$ denote the sum of all principal 
$2i$-minors of $M$.
If $n+m$ is odd  (resp. even), then we take the basic invariants $f_1,\dots,f_l$
(resp. $f_1,\dots,f_{l-1}$, and the pfaffian $\text{Pf}$). Let us prove that they form a
good generating system in $\textit{Inv}(\mathfrak{so}_{n+m},\ad)$  for $(\g,\g_0)$.

It easily follows from the block structure of $\tilde M$ that
$f_i\vert_{\ce\oplus\es}$ has a monomial entirely in $d_i$'s if and only if 
$i\le m$. Furthermore, if $i>m$, then one can always find a monomial in
$f_i\vert_{\ce\oplus\es}$ whose degree with respect to $d_i$'s equals $2m$.
Thus,
$\bideg f_i^\bullet =\begin{cases} (0,2i) & \text{ if $i\le m$},\\
                                                     (2i-2m,2m) & \text{ if $i>m$}.
\end{cases}$   
Likewise, for $n+m$ even, we have the pfaffian, 
and $\bideg\text{Pf}^\bullet=((n-m)/2,m)$.   
Actually, it is easily seen that 

$f_i^\bullet(\tilde M)=\begin{cases}
\text{the $i$-th elementary symmetric function in $d_1^2,\dots,d_m^2$}, & i\le m ; \\
\displaystyle\bigl( \prod_{i=1}^m d_i^2\bigr){\cdot}f_{i-m}(E), & i>m .
\end{cases}$
\\
and $\text{Pf}^\bullet(\tilde M)
= \displaystyle\bigl(\prod_{i=1}^m d_i\bigr)
{\cdot}\text{Pf}(E)$. Consequently, the
$f_i^\bullet\vert_{\ce\oplus\es}$, $i=1,\dots, \lfloor (n{+}m{-}1)/2\rfloor$,
 (together with  $\text{Pf}^\bullet\vert_{\ce\oplus\es}$ if
$n+m$ is even) are algebraically independent.
\end{proof}%
The following case is rather similar, although a bit more involved.

\begin{s}{Theorem}  \label{good:gl}
There is a good generating system for $(\g,\g_0)=
(\mathfrak{gl}_{n+m},\mathfrak{gl}_m\dotplus\gln)$.
\end{s}\begin{proof}
Here $l=n+m$ and $\rk(\g,\g_0)=\min\{n,m\}$. 
We use the natural matrix model for $(\g,\g_0)$:

$\g_0=\left\{\begin{pmatrix} M_1 & 0 \\ 0 & M_4 \end{pmatrix}\right\}$ and 
$\g_1=\left\{\begin{pmatrix} 0 & M_2 \\ M_3 & 0 \end{pmatrix}\right\}$,
\\[.6ex]
where $M_1$ (resp. $M_4$) is a  matrix of order $m$ (resp. $n$),
$M_2$ is a $m\times n$ matrix, and $M_3$ is a $n\times m$ matrix.
Assume below that $n\ge m$. Then 
$\es\simeq\mathfrak{gl}_{n-m}\dotplus\te_m$. 

Let us describe our choice of $\ce\subset\g_1$ and thereby of
$\es=\z_{\g_0}(\es)$. We take $M_2,M_3$ such that
$M_2=\begin{pmatrix} B & 0 \end{pmatrix}$ and
$M_3=\begin{pmatrix} -B \\ 0 \end{pmatrix}$, where $B$ is an 
arbitrary diagonal $m\times m$ matrix. 

Then taking the partition of $M$ into the nine submatrices
corresponding to the sizes  $m,m,n-m$, we obtain
\begin{equation}   \label{m_form}
\ce\oplus\es=\left\{\tilde M=\begin{pmatrix} A & B & 0 \\  -B & A & 0 \\ 0 & 0 & E 
\end{pmatrix}\right\} ,
\end{equation}
where $A$ and $B$ are diagonal matrices of order $m$ and
$E$ is an arbitrary matrix of order $n-m$. Let $a_1,\dots,a_m$ 
(resp. $b_1,\dots,b_m$) be the diagonal
entries of $A$ (resp. $B$). 

Let $f_i(M)$ denote the sum of all principal minors of order $i$ of a square matrix 
$M$. Let us prove that $f_1,\dots, f_{n+m}$ form a good generating system in 
$\textit{Inv}(\mathfrak{gl}_{n+m},\ad)$.  
It easily follows from 
Eq.~\eqref{m_form} that the restriction of $f_i$ 
to $\ce\oplus\es$ has a monomial entirely in $b_i$'s if and only if
$i$ is even and $i\le 2m$. If $i$ is odd and $i<2m$, then one can only find
a monomial whose all but one indeterminates are some
$b_i$'s. One other indeterminate is either
an $a_j$ (where $j$ depends on the $b_i$'s chosen) or an arbitrary diagonal entry of $E$.
Finally, if $i>2m$, then
one can always produce a monomial of $f_i\vert_{\ce\oplus\es}$ 
whose degree with respect to $b_i$'s equals $2m$.
Thus,
$\bideg f_i^\bullet =\begin{cases} (0,i) & \text{ if $i\le 2m$ and $i$ is even,}\\
                                                      (1,i-1) & \text{ if $i < 2m$ and $i$ is odd,}\\
                                                      (i-2m,2m) & \text{ if $i>2m$.}
\end{cases}$
\\[.7ex]
To describe these polynomials explicitly, we need some notation. Let $\sigma_i$
denote the $i$-th elementary symmetric function. Set 
$\mathcal{A}=\begin{pmatrix} A & 0 \\ 0 & A \end{pmatrix}$ and
$\mathcal{B}=\begin{pmatrix} 0 & B\\ -B & 0 \end{pmatrix}$.
Then looking at the principal minors of $\tilde M$ and their highest components with respect
to $\ce$ (i.e., $B$), one easily obtains

$\begin{cases}
f_{2i}^\bullet(\tilde M)=\sigma_i(b_1^2,\dots,b_m^2), & i\le m ; \\
f_{2i+1}^\bullet(\tilde M)=\tr(E)\sigma_i(b_1^2,\dots,b_m^2)+
\tr(\mathcal{B}^{2i}\mathcal{A}),  & i< m ;\\
f_j^\bullet(\tilde M)=f_{j-2m}(E)\sigma_m(b_1^2,\dots,b_m^2), & j> 2m.
\end{cases}$
\\[.6ex]
These formulae show that the polynomials $\{f_i^\bullet\}$ are algebraically independent.
\end{proof}%
Another example concerns an exceptional symmetric pair. This was obtained in collaboration
with O.\,Yakimova.

\begin{s}{Theorem}
For\/ $(\g,\g_0)=(\GR{F}{4},\GR{B}{4})$, there is a good generating system. 
The bi-degrees of the basic invariants in\/ $\bbk[\ka^*]^K$ are
$(0,2),\,(2,4),\,(4,4),\,(6,6)$.
\end{s}%
{\sl Sketch of the proof}. We explicitly construct a good generating system, using
{\sl ad hoc\/} arguments. Here $\dim\ce=1$, $\es\simeq \GR{B}{3}$, 
and we work with the restrictions of $G$-invariant functions to $\te=\ce\oplus\te(\es)$. 
The latter is a Cartan subalgebra of $\g$ and, by virtue of Chevalley's restriction theorem,
we actually deal with the Weyl group invariants on it.
The Weyl group of $\GR{F}{4}$,
$W(\GR{F}{4})$,  is a semi-direct product of the normal subgroup 
$W(\GR{D}{4})$ and $\mathcal{S}_3=W(\GR{A}{2})$. Here $W(\GR{D}{4})$ is generated
by the reflection with respect to the long roots of $\GR{F}{4}$ and $W(\GR{A}{2})$
is generated by the reflections corresponding to the short simple roots of $\GR{F}{4}$. 
Hence, to obtain $W(\GR{F}{4})$-invariants, one can take the invariants
of $W(\GR{D}{4})$ and then consider the $\mathcal S_3$-action on them.
We begin with a natural set of basic invariants of $W(\GR{D}{4})$.
The $\mathcal S_3$-action has a rather bulky expression with respect to this
set, but it is still a manageable 
task to write explicitly down the expressions for $W(\GR{F}{4})$-invariants
through the $W(\GR{D}{4})$-invariants. Then, playing around with these invariants,
we "correct" them on order to obtain a good generating system.

Here are the relevant data. 
We use the expressions for the simple roots of $\GR{F}{4}$
and their numbering from \cite{vo}; that is,
$\ap_1=\frac{1}{2}(\esi_1-\esi_2-\esi_3-\esi_4)$, $\ap_2=\esi_4$, $\ap_3=\esi_3-\esi_4$,
and $\ap_4=\esi_2-\esi_3$. Then $\Delta(\GR{D}{4})=\{\pm\esi_i\pm\esi_j\mid i< j\}$. 
The Satake diagram of our symmetric pair is:
\begin{picture}(78,18)(20,5)
\multiput(50,8)(20,0){3}{\circle*{6}}
\put(30,8){\circle{6}}
\put(33,8){\line(1,0){14}}
\put(73,8){\line(1,0){14}}
\multiput(52.5,7)(0,2){2}{\line(1,0){15}}
\put(55,5){$<$} 
\end{picture}.  (The white node represents $\ap_1$).
This shows that the simple roots of $\es$ are $\ap_2,\ap_3,\ap_4$ and
allows us to determine the splitting of $\te_{\BR}$.
Here $\ce_\BR=\BR\esi_1$ and $\te(\es)_\BR=
\BR\esi_2\oplus\BR\esi_3\oplus\BR\esi_4$. The basic invariants of $W(\GR{D}{4})$
are:
\begin{gather*}
 f_2=\esi_1^2+\esi_2^2+\esi_3^2+\esi_4^2 ,\\
 f'_4=\esi_1\esi_2\esi_3\esi_4, \\
 f_4=\sum_{i<j}\esi_i^2\esi_j^2, \\
 f_6=\sum_{i<j<k}\esi_i^2\esi_j^2\esi_k^2
\end{gather*}
The group $\mathcal S_3$ is realised as the group generated by the reflections $s_{\ap_1}$
and $s_{\ap_2}$. Using this, it straightforward to write down the $\mathcal S_3$-action
on $W(\GR{D}{4})$-invariants and to determine the $W(\GR{F}{4})$-invariants.
The key observation is that $f_2$ and $f_6-\frac{1}{6}f_2f_4$ are already 
$W(\GR{F}{4})$-invariants, and the plane $\textrm{Span}\{ f'_4, 4f_4-f_2^2\}$ 
affords the standard 
reflection representation of $\mathcal S_3$.
The basic invariants of $\GR{F}{4}$ have degrees $2,\,6,\,8,\,12$. Here are the expressions
of a good generating system $g_2,g_6,g_8,g_{12}$ via the $f_i$'s:
\begin{gather*} 
g_2= f_2, \\
g_6=f_6-\frac{1}{6}f_2f_4, \\
g_8={f'_4}^2+\frac{1}{12}f_4^2-\frac{1}{4}f_2f_6, \\
g_{12}=4{f'_4}^2f_4-\frac{3}{2}f_6^2-\frac{3}{2}{f'_4}^2f_2^2-\frac{1}{9}f_4^3
+\frac{1}{2}f_2f_4f_6 .
\end{gather*}
The highest components of these polynomials with the respect to $\ce$, i.e., 
with respect to $\esi_1$ are:
\begin{gather*}
g_2^\bullet= \esi_1^2, \\
g_6^\bullet= \esi_1^4(\esi_2^2+\esi_3^2+\esi_4^2), \\
g_8^\bullet=  \esi_1^4\bigl(\frac{1}{12}(\esi_2^2+\esi_3^2+\esi_4^2)^2-
\frac{1}{4}(\esi_2^2\esi_3^2+\esi_2^2\esi_4^2+\esi_3^2\esi_4^2)\bigr), \\
g_{12}^\bullet= \esi_1^6\bigl(-\frac{3}{2}\esi_2^2\esi_3^2\esi_4^2-
\frac{1}{9}(\esi_2^2+\esi_3^2+\esi_4^2)^3+
\frac{1}{2}(\esi_2^2+\esi_3^2+\esi_4^2)(\esi_2^2\esi_3^2+\esi_2^2\esi_4^2+
\esi_3^2\esi_4^2)\bigr).
\end{gather*}
It follows that these highest components are algebraically independent. \hfill $\square$
\\[.8ex]
It is likely that, for all symmetric pairs not mentioned in Remark~\ref{bad}, there is
a good generating system. However, this is not easy to prove, even for the other classical
series.


\section{$\N$-regular $\BZ_2$-gradings and their contractions}
\label{sect:N-reg}
\setcounter{equation}{0}

\noindent
A $\BZ_2$-grading (a symmetric pair) is said to be $\N$-{\it regular\/} if $\g_1$ 
contains a regular nilpotent element of $\g$.
By \cite{an}, 
a $\BZ_2$-grading is  $\N$-regular if and only if $\g_1$ contains a regular
semisimple element if and only if any nilpotent $G$-orbit in $\g$ meets $\g_1$.
 (This is no longer true for $\BZ_m$-gradings with $m>2$.)

Until the end of this section, we assume that our $\BZ_2$-grading is $\N$-regular.
Let $\ce\subset\g_1$ be a Cartan subspace. 

Set $Z_1=\ov{G{\cdot}\g_1}=\ov{G{\cdot}\ce}$. By \cite[Theorem\,4.7]{theta},
$Z_1$ is a normal complete intersection in $\g$ and the ideal of $Z_1$ in
$\bbk[\g]$ is generated by certain basic invariants. 
That is, 
there is a set of basic invariants
$f_1,\ldots,f_l$ such that $f_i\vert_{Z_1}\equiv 0$ for $i\ge k+1$
and $\bbk[Z_1]^G$ is freely generated by $f_i\vert_{Z_1}$ for $i\le k$.
Furthermore, since the restriction map $\bbk[\g]^G\to \bbk[\g_1]^{G_0}$ is 
onto~\cite[Theorem\,3.5]{theta},
$\bbk[\g_1]^{G_0}$ is freely generated by $\bar f_i=f_i\vert_{\g_1}$, 
$i\le k$, and $k=\rk(\g,\g_0)$. 
Thus, each $f_i$, $i=1,\dots,k$, has the bi-homogeneous component that does not depend on
$\g_0$, whereas $f_j$, $j=k{+}1,\dots,l$, does not have such a bi-homogeneous component.

If $x\in\g_1\cap\g_{reg}$, 
then $(\textsl{d}f_i)_x$, $i=1,\dots,l$, are linearly independent \cite{ko63}. 
Because  $f_j$ \ ($k{+}1{\le} j{\le} l$) does not have a  
component of degree 0 with respect to $\g_0$, it must have a
component of degree 1 with respect to $\g_0$. Otherwise we would have 
$(\textsl{d}f_j)_v=0$ for any $v\in\g_1$.
The linear component of $f_j$ with respect to $\g_0$ can be written as
$(x_0,x_1)\mapsto\langle x_0, F_j(x_1)\rangle$, 
where $0\ne F_j\in\Mor(\g_1,\g_0)$ and $\deg F_j=\deg f_j-1$. Since each bi-homogeneous
component of $f_j$ is $G_0$-invariant, $F_j$ must be $G_0$-equivariant, i.e.,
$F_j\in\Mor_{G_0}(\g_1,\g_0)$, $j=k{+}1,\dots,l$. 

As $\dim\g_1\md G_0=k$, the $\N$-regularity implies that 
$\dim\g_{1,x}=k$ and $\dim\g_{0,x}=l-k$ whenever $x\in\g_1\cap\g_{reg}$.
This also shows that $\dim\g_1-k=\dim\g_0-(l-k)$.
In view of $G_0$-equivariance, $F_j(x)\in \g_{0,x}$,
and the linear independence of the differentials 
$(\textsl{d}f_i)_x$ imply that $\{F_j(x)\}$  are linearly independent.
Hence  $\{F_j(x) \mid j=k{+}1,\dots,l\}$ is a basis for $\g_{0,x}$
for any  $x\in\g_1\cap\g_{reg}$.
Thus, we obtain the following presentation of the basic invariants 
$f_1,\dots,f_l$:
\begin{equation}  \label{presentation}
\left\{
\begin{array}{ll}
   f_i(x_0,x_1)=\bar f_i(x_1) + \text{(terms of higher degree w.r.t. $x_0$)}, \quad i\le k;\\
   f_j(x_0,x_1)=\langle x_0, F_j(x_1)\rangle  + 
\text{(terms of higher degree w.r.t. $x_0$)},  \quad j\ge k+1 \ .
\end{array}\right.
\end{equation}
Set $\widehat F_j(x_0,x_1)=\langle x_0, F_j(x_1)\rangle$.

\begin{s}{Theorem}   \label{n-reg}
Let\/ $\ka=\g_0\ltimes\g_1$ be the $\BZ_2$-contraction of an $\N$-regular $\BZ_2$-grading
of rank $k$. Then, using the above notation,
\begin{itemize}
\item[\sf (i)] \ $\bbk[\ka^*]^{K^u}$ is a polynomial algebra that is freely
generated by the coordinates on $\g_1^*$ and  $\widehat F_j$, $j=k+1,\dots,l$.
\item[\sf (ii)] \ $\bbk[\ka^*]^K$ is the polynomial algebra 
that is freely generated by $\bar f_1,\dots,\bar f_k,\widehat F_{k+1},\dots,
\widehat F_l$.
\end{itemize}
\end{s}\begin{proof} 
It follows from Eq.~\eqref{presentation} that
$f_i^\bullet=\bar f_i$ for $i\le k$ and $f_j^\bullet=\widehat F_j$ for $j\le k+1$.
Regarding all these functions as functions on $\ka^*$, we obtain by virtue of 
Proposition~\ref{z2-contra} that 
$\bar f_1,\dots,\bar f_k,\widehat F_{k+1},\dots,\widehat F_l$ belong to 
$\bbk[\ka^*]^K$.

The proof below is quite similar to that of Theorem~6.2 in \cite{p05}.
To prove part (i), we use Igusa's lemma (see \cite[Theorem\,4.12]{t55}
or \cite[Lemma\,6.1]{p05}) and 
properties of the covariants $F_i$ in Eq.~\eqref{presentation}. Part (ii)
is then an obvious consequence of (i).

For (i): Let $\Omega\subset \g_1^*$ be the open subset of $G_0$-regular elements.
As follows from \cite{kr}, $\Omega$ is big and any $\zeta\in\Omega$ is also regular as 
element of $\g$. The functions indicated in (i) are clearly $K^u$-invariant and are 
algebraically independent (consider their differentials at some $\xi\in\Omega$).
Hence we obtain the dominant mapping
\[
  \psi:\ka^* \to \g_1^* \times \bbk^{l-k},
\]
defined by $\psi(\xi_0,\xi_1)=(\xi_1, \widehat F_{k+1}(\xi_0,\xi_1),\dots,\widehat F_l(\xi_0,\xi_1))$.
Since the vectors $F_j(\zeta)$, $j=k{+}1,\dots,l$, form a 
basis of $\g_{0,\zeta}$ for any $\zeta\in\Omega$, we see that
$\Omega\times \bbk^{l-k} \subset \Ima\psi$, i.e., $\Ima\psi$ contains
a big open subset of $\g_1^*\times\bbk^{l-k}$.
If $(\zeta,z_{k+1},\dots,z_l)\subset \Omega\times \bbk^{l-k}$, then 
\[
   \psi^{-1}(\zeta,z_{k+1},\dots,z_l)=\{ (\xi_0,\zeta)\mid \langle \xi_0, F_j(\zeta)\rangle
=z_j, \ j\ge k+1 \} .
\]
It is a $K^u$-stable affine subspace of $\ka^*$ of dimension $\dim\g_0-(l-k)$. 
On the other hand,
if $(\xi_0,\zeta)\in \psi^{-1}(\zeta,z_{k+1},\dots,z_l)$, then 
\[
 K^u{\cdot}(\xi_0,\zeta)=(1\ltimes\g_1){\cdot}(\xi_0,\zeta)=
\{(\xi_0 + x_1\ast\zeta,\zeta) \mid x_1\in \g_1 \} \ .
\]
Upon the identification of $\g_1$ and $\g_1^*$, we have $\g_1\ast\zeta= [\g_1,\zeta]$.
Hence 
\[
\dim K^u{\cdot}(\xi_0,\zeta)=\dim (\g_1{\ast}\zeta)=
\dim\g_1-k=\dim\g_0-(l-k) \ .
\]
Since the orbits of  unipotent groups on affine varieties are closed 
\cite[Theorem\,2]{ros61}
and isomorphic to affine spaces, we conclude that 
$\psi^{-1}(\zeta,z_{k+1},\dots,z_l)=K^u(\xi_0,\zeta)$, i.e., almost all fibres of
$\psi$ are precisely $K^u$-orbits.

Hence all the assumptions of Igusa's lemma are satisfied, and part (i) follows.

A direct proof for part (ii) (without using (i)) is as follows.
The $K$-invariants $\bar f_1,\dots,\bar f_k,\widehat F_{k+1},\dots,\widehat F_l$
are $\BZ_2$-degenerations of $f_1,\dots,f_l$, hence they have the same degrees.
It is also easily seen that these $K$-invariants are algebraically independent.
Next, we know that $\ind\ka=\ind\g$. Therefore, Theorem~\ref{ner-vo}(ii)
applies in this situation. 
\end{proof}%
{\bf Remark.} If $\rk(\g,\g_0)=l=\rk\g$, then the above theorem merely says that
$\bbk[\ka^*]^{K^u}\simeq \bbk[\g_1]$ and $\bbk[\ka^*]^K\simeq \bbk[\g_1]^{G_0}$.
This was already observed, in a more general context, in \cite[Theorem\,6.4]{p05}.
So, the novelty of Theorem~\ref{n-reg} concerns the case in which 
$\rk(\g,\g_0) < \rk\g$.

\begin{s}{Theorem}  \label{n-eq}
If $\ka=\g_0\ltimes\g_1$ is the $\BZ_2$-contraction of an $\N$-regular $\BZ_2$-grading, 
then $\pi_{\ka^*}: \ka^* \to \ka^*\md K$ is equidimensional.
\end{s}\begin{proof}
Keep the notation of the previous proof.
If $\rk(\g,\g_0)=\rk\g$, then the isomorphism $\bbk[\ka^*]^K\simeq \bbk[\g_1]^{G_0}$
shows that $\gN(\ka^*)\simeq \gN(\g_1)\times\g_0$. Hence the assertion.

Assume therefore that $k=\rk(\g,\g_0) <\rk\g=l$ and hence there are non-trivial
$K$-invariants of the form $\widehat F_j$. Roughly speaking,
the covariants $F_j\in\Mor_{G_0}(\g_1,\g_0)$, $j=k{+}1,\dots,l$, determine a stratification
of the null-cone ${\goth N}(\g_1)$, and the assertion
is equivalent to certain property of this stratification. 
Unfortunately, we can only verify that property using a case-by-case
argument. (This is very similar to our proofs for $(\ka,\ad)$ 
in \cite[Sect.\,9]{p05}.) 

Set $\gN(\ka^*)=\pi_{\ka^*}^{-1}(\pi_{\ka^*}(0))$. Then $\codim \gN(\ka^*)\le l$ and the 
equidimensionality of $\pi_{\ka^*}$ precisely means that $\codim \gN(\ka^*)= l$. 
If $(\ap,\xi)\in\gN(\ka^*)$, then
the inclusion $\bbk[\g_1]^{G_0}\hookrightarrow\bbk[\ka^*]^K$ shows that $\xi\in\gN(\g_1)$.
Hence we obtain the surjective projection 
$p: \gN(\ka^*)\to \gN(\g_1)$, $(\ap,\xi)\mapsto \xi$.
Let $J$ denote the (finite) set of $G_0$-orbits in $\gN(\g_1)$. Then 
$\gN(\ka^*)=\sqcup_{\co\in J}p^{-1}(\co)$ and the irreducible components
are contained among the sets $\ov{p^{-1}(\co)}$. Hence the assertion is equivalent to the 
condition that $\dim p^{-1}(\co)\le \dim\ka^*-l$ for all $\co\in J$.
Since $\dim p^{-1}(\xi)=\dim\g_0- \dim \text{span}\{F_{k+1}(\xi),\dots,F_{l}(\xi)\}$
for $\xi\in \gN(\g_1)$, 
%
our condition readily translates  as follows:
For any  $\xi\in\gN(\g_1)$, we should have
\begin{equation}  \label{uslovie}
   l \le \dim \g_{1,\xi}+\dim \text{span}\{F_{k+1}(\xi),\dots,F_{l}(\xi)\} \ .
\end{equation}
Recall that $\dim\g_{1,\xi}\ge k=\rk(\g,\g_0)$ and $\dim\g_{1,\xi}=k$ if and only if
$\xi\in\Omega$.

The list of $\N$-regular symmetric pairs such that $\g$ is simple 
and $\rk\g> \rk(\g,\g_0)$ is given below:

$(\mathfrak{sl}_{n+m},\sln\dotplus\mathfrak{sl}_m\dotplus\bbk)$ \ with $|n-m|\le 1$;

$(\mathfrak{so}_{2n+2},\son\dotplus\mathfrak{so}_{n+2})$;

$(\GR{E}{6}, \mathfrak{sl}_6\dotplus\tri)$.
\\[.6ex]
In the second case, $l=n+1$ and $k=n$. Here there is only one covariant, $F_l$, 
and the equidimensionality is obvious.

In the third case, $l-k=6-4=2$, and there are two covariants $F_5,F_6$. 
Here Eq.~\eqref{uslovie} essentially
means that  if $\xi\in \gN(\g_1)$ and 
$\codim_{\gN(\g_1)}G_0{\cdot}\xi=1$, then at least one of the 
covariants $F_5, F_6$ does not vanish at $\xi$.
To this end, we notice that  $G{\cdot}\xi$ is the subregular nilpotent orbit, $\co_{sub}$, 
For $\xi\in\co_{sub}$, we have 
$\dim\text{span}\{(\textsl{d}f_1)_\xi,\dots,(\textsl{d}f_l)_\xi\}=l-1$ 
\cite{bram}.
For our "adapted" choice of basic invariants $f_1,\dots,f_l$, as above,
we have $(\textsl{d}f_j)_\xi=F_j(\xi)$ for $j\ge k+1$ and $\xi\in\g_1$. 
Hence two covariants $F_j$ cannot vanish on $\co_{sub}\cap \g_1$, which is exactly
what we need.

For the first case, Eq.~\eqref{uslovie} will be verified in 
Example~\ref{primery-reg} below.
%
\end{proof}%
The following is a standard consequence of 
Theorem~\ref{n-reg}(ii) and Theorem~\ref{n-eq}.

\begin{s}{Corollary}
Let\/ $\ka=\g_0\ltimes\g_1$ be the contraction of an $\N$-regular $\BZ_2$-grading. 
Let ${\goth U}(\ka)$ denote the enveloping algebra of $\ka$ and ${\goth Z}(\ka)$
the centre of\/ ${\goth U}(\ka)$. Then ${\goth Z}(\ka)$ is a polynomial algebra 
and\/ ${\goth U}(\ka)$ is a free module over ${\goth Z}(\ka)$.
\end{s}%
\vskip-1ex
\begin{rem}{Example}   \label{primery-reg}
To simplify exposition, we work with $\gln$ in place of $\sln$. 
Let $\g=\mathfrak{gl}_{2n}$ and 
$\g_0=\gln\dotplus\gln=\begin{pmatrix} \ast & 0 \\ 0 & \ast \end{pmatrix}$. Then
$\g_1=\begin{pmatrix}  0 & \ast \\ \ast & 0 \end{pmatrix}$.
As always, we identify $\g_i$ and $\g_i^*$, $i=0,1$.
Here $\rk\g=2n$ and $\rk(\g,\g_0)=n$.
Set $\xi_0=
\begin{pmatrix} M & 0 \\ 0 & N \end{pmatrix}\in\g_0$  and 
$\xi_1=
\begin{pmatrix} 0 & A \\ B & 0 \end{pmatrix}\in\g_1$.
To obtain a regular nilpotent element in $\g_1$, one may take $B=I_n$ and $A$ 
to be any nilpotent $n\times n$ matrix such that $A^{n-1}\ne 0$.
The algebra $\bbk[\g_1]^{G_0}$ is freely generated by the polynomials
$f_i(\xi_1)=\tr((\xi_1)^{2i})=\tr((AB)^i+(BA)^i)$, $i=1,2,\dots,n$.
These polynomials are naturally regarded as polynomials on the whole of $\ka^*$.
Define the covariants $F_i:\g_1\to \g_0$ by the formula
\[
   F_i(\xi_1)=\begin{pmatrix} 0 & A \\ B & 0 \end{pmatrix}^{2i-2}
=\begin{pmatrix} (AB)^{i-1} & 0 \\ 0 & (BA)^{i-1} \end{pmatrix} .
\]
Obviously, $F_i(\xi_1)$ commutes with $\xi_1$, i.e., $F_i(\xi_1)\in \g_{0,\xi_1}$.
Therefore $\widehat F_i(\xi_0,\xi_1):=\langle \xi_0, F_i(\xi_1)\rangle=
\tr(\xi_0 (\xi_1)^{2i})$ is a $K$-invariant polynomial on $\ka^*$.
If $\xi_1\in\g_1$ is a regular nilpotent element, then
$\xi_1^0,\xi_1^2,\dots,\xi_1^{2n-2}$ are linearly independent. Hence $F_1,\dots,F_n$ form a 
basis of the $\bbk[\g_1]^{G_0}$-module $\Mor_{G_0}(\g_1,\g_0)$.
It follows that $\bbk[\ka^*]^K$ is freely generated by the polynomials

$f_i(\xi_0,\xi_1)=\tr((AB)^i+(BA)^i)$, $i=1,2,\dots,n$, \quad  and

$\widehat F_i(\xi_0,\xi_1)=\tr(M(AB)^i+N(BA)^i)$, $i=1,\dots,n$.
\\[.6ex]
In this case, Eq.~\eqref{uslovie} for $\xi\in\gN(\g_1)$ reads
\[
    \dim\g_{1,\xi} +\dim\text{span}\{\xi^{2i} \mid i=0,1,\dots,n-1\}-2n \ge 0 \ .
\]
Eliminating $i=0$ and taking into account that here $\dim\g_{1,\xi}=
\frac{1}{2}\dim\g_\xi$, we rewrite it as
\[
  \frac{1}{2}\dim\g_\xi+\dim\text{span}\{\xi^{2i} \mid i=1,\dots,n-1\}-2n+1 \ge 0 \ .
\]
Let $(\eta_1,\eta_2,\ldots)$ be the  partition of $2n$ corresponding
to $\xi$. Then $\xi^{2i}\ne 0$ if and only if $2i\le \eta_1-1$.
Write $(\hat\eta_1,\hat\eta_2,\ldots,\hat\eta_s)$ for the dual partition. 
This means in particular that $s=\eta_1$. 
It is well-known that $\dim\g_\xi=\sum_{i=1}^s\hat\eta_i^2$. Hence the left-hand side
equals
\begin{multline*}
\frac{1}{2}\sum_{i=1}^s\hat\eta_i^2+\left\lfloor\frac{\eta_1-1}{2}\right\rfloor-2n+1=
\frac{1}{2}\sum_{i=1}^s\hat\eta_i^2 +\left\lfloor\frac{s-1}{2}\right\rfloor
-(\sum_{i=1}^s \hat\eta_i) + 1 = 
\\
\frac{1}{2}\bigl(\sum_{i=1}^s(\hat\eta_i-1)^2 -s+2 
+2\left\lfloor\frac{s-1}{2}\right\rfloor\bigr)= 
\frac{1}{2}\sum_{i=1}^s(\hat\eta_i-1)^2 +
\bigl(\left\lfloor\frac{s+1}{2}\right\rfloor-\frac{s}{2}\bigr)
 \ ,
\end{multline*}
which is nonnegative, as required.
\\[.6ex]
The case of $\g=\mathfrak{gl}_{2n+1}$ and $\g_0=\gln\dotplus \mathfrak{gl}_{n+1}$
is quite similar and left to the reader.
\end{rem}%
\vskip-1ex

\begin{rem}{Remarks}  \label{radical_ideal}
1. \ Using a more involved analysis, we can prove that, for all $\N$-regular $\BZ_2$-gradings,
the ideal generated by the basic $K$-invariants $\bar f_1,\dots,\bar f_k,\widehat F_{k+1},\dots,
\widehat F_l$ is equal to its radical. To this end, it suffices to demonstrate that
each irreducible component of $\gN(\ka^*)$ contains a $K$-regular point.
\\[.6ex]
2. \ The null-fibre $\gN(\ka^*)$ is often reducible. The projection 
$p:\gN(\ka^*)\to \gN(\g_1)$ considered in Theorem~\ref{n-eq} shows that
$\#\textrm{Irr}(\gN(\ka^*))\ge \#\textrm{Irr}(\gN(\g_1))$, where $\#\textrm{Irr}(\cdot)$
refers to the number of irreducible components. The numbers $\#\textrm{Irr}(\gN(\g_1))$
are found by Sekiguchi for all symmetric pairs~\cite[Theorem\,1]{jiro}. It may happen that
$\#\textrm{Irr}(\gN(\ka^*)) > \#\textrm{Irr}(\gN(\g_1))$. For instance, if 
$\g=\mathfrak{gl}_{2n+1}$ and $\g_0=\gln\dotplus \mathfrak{gl}_{n+1}$, then
$\gN(\g_1)$ is irreducible, while our computation shows that 
$\#\textrm{Irr}(\gN(\ka^*))=2$. The additional irreducible component appears as
the closure of $p^{-1}(\co_{sub}\cap\g_1)$.
\end{rem}%
The covariants $F_{k+1},\dots,F_l$ have another natural description.
Let $\Mor(\g_1,\g_0)$ (resp. $\Mor(\g_1,\g_1)$) be the set of {\sl all\/} polynomial
morphisms $\g_1\to \g_0$ (resp. $\g_1\to\g_1$).
These are free $\bbk[\g_1]$-modules
of rank $\dim\g_0$ and $\dim\g_1$, respectively. 
Consider the homomorphism
${\hat\phi}:\Mor(\g_1,\g_0){\to} \Mor(\g_1,\g_1)$
defined by $\hat\phi(F)(x_1)=[F(x_1),x_1]$. Then $\ker\hat\phi=\{F\mid F(x_1)\in\g_{0,x_1}\}$.
Notice that $F$ is not supposed to be $G_0$-equivariant. The homomorphism $\hat\phi$ can be defined for any symmetric pair. But the following is only true in the $\N$-regular case
 (cf. \cite[Theorem\,8.6]{p05}).

\begin{s}{Theorem}  \label{kernel}
Let $(\g,\g_0)$ be an $\N$-regular symmetric pair. Then
$\ker\hat\phi$ is a free $\bbk[\g_1]$-module, and the $G_0$-equivariant morphisms
$F_{k+1},\dots,F_l$ form a basis of\/ $\ker\hat\phi$.
\end{s}\begin{proof}
Clearly, $\ker\hat\phi$ is a torsion-free $\bbk[\g_1]$-module and 
its rank, $\rk(\ker\hat\phi)$, is well-defined. By definition, 
$\rk(\ker\hat\phi):=\dim(\ker\hat\phi\otimes_{\bkk[\g_1]} \bbk(\g_1))$.  
In the coordinate form, 
$\hat\phi$ is represented via a $\dim\g_0\times \dim\g_1$-matrix 
with entries in $\bbk[\g_1]$, and $\rk(\hat\phi)$ is the rank of this matrix.
Then $\rk(\ker\hat\psi)=\dim\g_0-\rk(\hat\psi)$. Because
\[
\rk\hat\phi=\max_{x\in\g_1}\dim G_0{\cdot}x=\dim\g_1-k=\dim\g_0-(l-k) \ ,
\]
we have $\rk(\ker\hat\psi)=l-k$.
Recall that  $F_{k+1}(\xi),\dots,F_l(\xi)$ are linearly independent over $\bbk$ for any
$\xi\in\Omega$,  hence $F_{k+1},\dots,F_l$ are linearly
independent over $\bbk[\g_1]$. 
As was noticed before, $F_{k+1},\dots,F_l\in\ker\hat\phi$.
Hence $F_{k+1},\dots,F_l$ generate $\ker\hat\phi\otimes_{\bkk[\g_1]} \bbk(\g_1)$.
That is,  for any $F\in \ker\hat\psi$ there exist $\hat p,p_{k+1},\ldots,p_{l}\in \bbk[\g_1]$ such that 
\[
     \hat pF=\sum_{i\ge k+1} p_i F_i \ . 
\]
Assume $\hat p\not\in \bbk^*$. Let $p$ be a prime factor
of $\hat p$ and $D$  the divisor of zeros of $p$.
Then $\sum_i p_i(v) F_i(v)=0$ for any $v\in D$. Since $\Omega\subset\g_1$ is big,
$\Omega\cap D$ is dense in $D$. Because 
$\{F_i(v)\}$ are linearly independent for any $v\in\Omega$, we obtain
$p_i\vert_D\equiv 0$. Hence $p_i/p\in \bbk[\g]$ for each $i$, and we are
done.
\end{proof}%
Note that $\ker\hat\phi$ cannot be generated by 
$G_0$-equivariant morphisms, unless $(\g,\g_0)$ is $\N$-re\-gu\-lar. The reason is that in general $\rk(\ker\hat\phi)=\dim\es$,
whereas one can show that the set of $G_0$-equivariant morphisms
in $\ker\hat\phi$ has the rank $\dim\z(\es)$ as the $\bbk[\g_1]^{G_0}$-module.
It remains to observe that $\N$-regularity precisely means that $\es$ is
commutative, i.e., $\es=\z(\es)$.


\section{Tables}
\label{sect:tables}
\setcounter{equation}{0}

\noindent
In this section, we gather the available information about the structure of
algebras $\bbk[\ka]^K$ and $\bbk[\ka^*]^K$, where $\ka=\g_0\ltimes\g_1$. 
The case of the adjoint representation is fully
covered by results of \cite{p05}. In particular, if 
$\rk\g=\rk\g_0$, i.e., the involution $\sigma$ is inner,
then $\bbk[\ka]^K\simeq \bbk[\g_0]^{G_0}$. Therefore we do not always write explicitly down
the respective bi-degrees. For $\bbk[\ka^*]^K$, the answer is known for the symmetric pairs 
considered in Sections~\ref{good_systems} and \ref{sect:N-reg}. In the maximal rank case,
we have $\bbk[\ka^*]^K\simeq \bbk[\g_1]^{G_0}\simeq \bbk[\g]^G$, and in two such cases
we omit indication of the degrees.
For all other pairs
not mentioned in Remark~\ref{bad}, we have precise suggestions for the degrees. 
These conjectural degrees are displayed in 
{\color{green}{\it italic}}. For the four cases mentioned in 
Remark~\ref{bad}, we put the question mark, if there is no suggestion for the corresponding 
degree, see Table~\ref{3}.

The last column contains a comment on the pair in question: "{\it max}" means the 
maximal rank case; $(N_i)$ means that $\g_i$ contains 
a regular nilpotent element of $\g$, $i=0,1$. [Hence $(N_1)=\N$-regular.]

Recall that $\rk\g=\rk(\g,\g_0)+ \rk\es$, $\dim\g\md G_0=\rk(\g,\g_0)$, and 
$\bbk[\g_1]^{G_0}$ is embedded in $\bbk[\ka^*]^K$. Hence we always have $\rk(\g,\g_0)$ 
basic invariants whose $\g_0$-degree  equals $0$.
There is an {\sl a posteriori\/} observation related to the $\g_0$-degrees 
of the remaining basic invariants in $\bbk[\ka^*]^K$. Namely, they are equal to the basic
degrees of $\es$ in all cases, where the algebra $\bbk[\ka^*]^K$ is known. 
For instance, look at the first pair in Table~\ref{1}. Here $\es\simeq\mathfrak{gl}_{n-m}
\times \te_m$ and the nonzero $\g_0$-degrees are $1$ ($m$ times), $1,2,3,\dots,n-m$.

\noindent  
\begin{table}[htb]
\caption{Classical Lie algebras (inner involutions)} \label{1}
\begin{tabular}{|c|cl|c|}   \hline
$(G,G_0)$  &  Bi-degrees for $(\ka,\ad)$ & Bi-degrees for $(\ka,\ads)$ 
&  \\ \hline
$(GL_{n+m}, GL_{n}\times GL_m)$ & degrees of $\g_0$ 
& (0,2),\,(0,4),\ldots,(0,2m) \rule{0pt}{2.5ex} &  $(N_1)$ if \\
n$\geqslant$m & & (1,0),\,(1,2),\,(1,4),\ldots,(1,2m) \rule{0pt}{2ex} & 
n-m$\leqslant$1  \\
 & & { (2,2m),\,(3,2m),\,\ldots,(n-m,2m)} \rule{0pt}{2ex} &  
\\
$(Sp_{2n+2m}, Sp_{2n}\times Sp_{2m})$ & degrees of $\g_0$ 
& (0,2),\,(0,4),\ldots,(0,2m) \rule{0pt}{2.5ex} &   \\
n$\geqslant$m & & {\color{green}{\it (2,2m),\,(2,2m+2),\,\ldots,(2,4m-2)}} 
\rule{0pt}{2.5ex} &   \\
   & & {\color{green}{\it (2,4m),\,(4,4m),\,\ldots,(2n-2m,4m)}} \rule{0pt}{2ex} &   
\\
$(SO_{2n+1}, SO_{n+1}\times SO_{n})$ & degrees of $\g_0$ 
& (0,2),\,(0,4),\ldots,(0,2n) \rule{0pt}{2.5ex} &  {\it max}  
\\
$(SO_{n+m}, SO_{n}\times SO_{m})$ & degrees of $\g_0$ 
& (0,2),\,(0,4),\ldots,(0,2m) \rule{0pt}{2.5ex} &   \\
n$>$m+2, n$+$m \ is odd & & 
{(2,2m),\,(4,2m),\,\ldots,(n-m-1,2m)} \rule{0pt}{2ex} &   
\\
$(SO_{4n}, GL_{2n})$ & degrees of $\g_0$ 
& (0,2),\,(0,4),\ldots,(0,2n) \rule{0pt}{2.5ex} &   \\
 & & 
{\color{green}{\it (2,2n-2),\,(2,2n),\,\ldots,(2,4n-4)}} \rule{0pt}{2.5ex} &  
\\
$(SO_{4n+2}, GL_{2n+1})$ & degrees of $\g_0$ 
& (0,2),\,(0,4),\ldots,(0,2n),\,(1,2n) \rule{0pt}{2.5ex} &   \\
 & & 
{\color{green}{\it (2,2n),\,(2,2n+2),\,\ldots,(2,4n-2)}} \rule{0pt}{2.5ex} &  
 \\
\hline
\end{tabular}
\end{table}

\noindent \begin{table}[htb]
\caption{Classical Lie algebras (outer involutions)}
\begin{tabular}{|c|ll|c|}   \hline
$(G,G_0)$  &  Bi-degrees for $(\ka,\ad)$ & Bi-degrees for $(\ka,\ads)$ 
&  \\ \hline
$(SL_{2n}, SO_{2n})$ &
(2,0),\,(4,0),\ldots,(2n-2,0),\,(n,0) & (0,2),\,(0,3),\ldots,(0,2n) & {\it max} \\
   &  (2,1),\,(4,1),\ldots,(2n-2,1) & &
\\
$(SL_{2n+1}, SO_{2n+1})$ &  
(2,0),\,(4,0),\ldots,(2n-2,0),\,(2n,0) 
& (0,2),\,(0,3),\ldots,(0,2n+1) \rule{0pt}{2.5ex} & {\it max} \\
   &  (2,1),\,(4,1),\ldots,(2n-2,1),\,(2n,1) &  &  $\&\, (N_0)$
\\
$(SL_{2n}, Sp_{2n})$ & 
(2,0),\,(4,0),\ldots,(2n-2,0),\,(2n,0) 
& (0,2),\,(0,3),\ldots,(0,n) \rule{0pt}{2.5ex} & $(N_0)$ \\
   &  (2,1),\,(4,1),\ldots,(2n-2,1) & {\color{green}{\it (2,n-1),\,(2,n),\ldots,(2,2n-2)}} & 
 \\
$(\GR{D}{n+m+1}, \GR{B}{n}{\times} \GR{B}{m})$ & 
(2,0),\,(4,0),\ldots,(2n-2,0),\,(2n,0) 
& (0,2),\,(0,4),\ldots,(0,4m+2) \rule{0pt}{2.5ex} & $(N_0)$ if\\
 n$\geqslant$m+2 & (2,0),\,(4,0),\ldots,(2m-2,0),\,(2m,0) 
   & { (2,4m+2),\,(4,4m+2),\ldots,} & m$=$0\\
   &  (n+m,1) & { (2n-2m-2,4m+2),\,(n-m,2m+1)} & 
\\ 
$(\GR{D}{2n}, \GR{B}{n}{\times} \GR{B}{n-1})$ & 
(2,0),\,(4,0),\ldots,(2n-2,0),\,(2n,0) 
& (0,2),\,(0,4),\ldots,(0,4n-2) \rule{0pt}{2.5ex} & $(N_1)$ \\
 (n=m+1) & (2,0),\,(4,0),\ldots,(2m-2,0),\,(2m,0) 
   & (1,2n-1) &  \\  &  (2n-1,1) &  & 
\\ 
$(\GR{D}{2n+1}, \GR{B}{n}{\times} \GR{B}{n})$ & 
(2,0),\,(4,0),\ldots,(2n-2,0),\,(2n,0) 
& (0,2),\,(0,4),\ldots,(0,4n),\,(0,2n+1) \rule{0pt}{2.5ex} & {\it max} \\
(n=m)  & (2,0),\,(4,0),\ldots,(2n-2,0),\,(2n,0) 
   &  &  \\  &  (2n,1) &  & \\
\hline
\end{tabular}
\end{table}

\vskip1ex


\noindent  \begin{table}[htb]   
\caption{Exceptional Lie algebras}  \label{3}
\begin{tabular}{|c|ll|c|}   \hline
$(G,G_0)$  &  Bi-degrees for $(\ka,\ad)$ & Bi-degrees for $(\ka,\ads)$ 
&  \\ \hline
$(\GR{F}{4}, \GR{B}{4})$ &  (2,0),\,(4,0),\,(6,0),\,(8,0) 
& (0,2),\,{(2,4),\,(4,4),\,(6,6)} \rule{0pt}{2.5ex} &
\\
$(\GR{F}{4}, \GR{C}{3}{\times}\GR{A}{1})$ &  (2,0),\,(2,0),\,(4,0),\,(6,0) & 
(0,2),\,(0,6),\,(0,8),\,(0,12)  &  {\it max} \rule{0pt}{2.5ex} 
\\
$(\GR{E}{6} , \GR{C}{4})$ & 
(2,0),\,(4,0),\,(4,1),\,(6,0),\,(8,0),\,(8,1) 
& (0,2),\,(0,5),\,(0,6),\,(0,8),\,(0,9),\,(0,12) &  {\it max} \rule{0pt}{2.5ex}
\\
$(\GR{E}{6} , \GR{F}{4})$ & (2,0),\,(4,1),\,(6,0),\,(8,0),\,(8,1),\,(12,0)
& (0,2),\,(0,3),\,{\color{green}{\it (2,?),\,(4,?),\,(4,?),\,(6,?)}} & $(N_0)$ 
\rule{0pt}{2.5ex}
\\
$(\GR{E}{6} , \GR{A}{5}{\times}\GR{A}{1})$   & 
(2,0),\,(2,0),\,(3,0),\,(4,0),\,(5,0),\,(6,0)  & (0,2),\,(0,6),\,(0,8),\,(0,12),\,(1,4),\,(1,8) & 
$(N_1)$ \rule{0pt}{2.5ex}
\\
$(\GR{E}{6} , \GR{D}{5}{\times}\GR{T}{1})$   & 
(1,0),\,(2,0),\,(4,0),\,(6,0),\,(8,0),\,(5,0)  & 
(0,2),\,(0,4),\,{\color{green}{\it (1,?),\,(2,?),\,(3,?),\,(4,?)}} & 
\rule{0pt}{2.5ex}
\\
$(\GR{E}{7} , \GR{D}{6}{\times}\GR{A}{1})$ & 
(2,0),\,(4,0),\,(6,0),\,(8,0),\,(10,0),\,(6,0) & (0,2),\,(0,6),\,(0,8),\,(0,12), 
 & \rule{0pt}{2.5ex}
\\  
& (2,0)  &  {\color{green}{\it (2,8),\,(2,12),\,(2,16)}}  & 
\\ 
$(\GR{E}{7} , \GR{E}{6}{\times}\GR{T}{1})$ & 
(1,0),\,(2,0),\,(5,0),\,(6,0),\,(8,0),\,(9,0) & (0,2),\,(0,4),\,(0,6), 
 & \rule{0pt}{2.5ex}
\\  
& (12,0)  &  {\color{green}{\it (2,?),\,(4,?),\,(4,?),\,(6,?)}}  & 
\\ 
$(\GR{E}{7} , \GR{A}{7})$ & 
degrees of $\GR{A}{7}$
& degrees of $\GR{E}{7}$ &  {\it max} \rule{0pt}{2.5ex}
\\
$(\GR{E}{8} , \GR{E}{7}{\times}\GR{A}{1})$ & 
(2,0),\,(6,0),\,(8,0),\,(10,0),\,(12,0) & 
(0,2),\,(0,6),\,(0,8),\,(0,12), 
 & \rule{0pt}{2.5ex} \\  
& (14,0),\,(18,0),\,(2,0)  &  {\color{green}{\it (2,?),\,(4,?),\,(4,?),\,(6,?)}}  & 
\\
$(\GR{E}{8} , \GR{D}{8})$ & 
degrees of $\GR{D}{8}$
& degrees of $\GR{E}{8}$ &  {\it max} \rule{0pt}{2.5ex}
\\
\hline
\end{tabular}
\end{table}

\end{document}